\documentclass[11pt]{article}
\usepackage[T1]{fontenc}
\usepackage[utf8]{inputenc}
\usepackage{fullpage,times}
\usepackage{microtype,needspace,enumitem}
\usepackage{amsmath,amssymb,amsthm,booktabs,array}
\usepackage[hyphens]{url}
\usepackage[hidelinks,bookmarksnumbered]{hyperref}
\hypersetup{pdftitle={The random k-SAT threshold is at most the survey propagation prediction},
  pdfauthor={Dimitris Achlioptas and Fedor Vorobyev},
  pdfsubject={Convergence of survey propagation and certified enclosures of the SP prediction}}
\newtheorem{theorem}{Theorem}
\newtheorem{proposition}[theorem]{Proposition}
\newtheorem{corollary}[theorem]{Corollary}
\newtheorem{lemma}[theorem]{Lemma}
\newcommand{\E}{\mathbb E}
\newcommand{\Prb}{\mathbb P}

\newcommand{\Pois}{\operatorname{Po}}
\newcommand{\Law}{\operatorname{Law}}

\newcommand{\PHI}{\boldsymbol{\Phi}}
\newcommand{\PHIaux}{\widehat{\PHI}}
\newcommand{\rk}{r_k}
\newcommand{\rstar}{r_{\mathrm{1RSB}}(k)}
\newcommand{\surv}{\mathfrak u}
\newcommand{\res}{\operatorname{res}}
\newcommand{\Lfun}{\mathcal L} 
\title{The random $k$-SAT threshold is at most \\ the survey propagation prediction}
\author{Dimitris Achlioptas\\University of Athens
  \and Fedor Vorobyev\\Independent researcher}
\date{}

\begin{document}
\maketitle
\begin{abstract}
We prove that the random $3$-SAT threshold is at most $4.26678429403$. More generally, we prove that for every $k\ge3$, the random $k$-SAT threshold is at most the value predicted by Survey Propagation. For all $k\ge 3$ and $b \ge 2$, we  show how to compute an interval of width $10^{-b}$ containing this value in $O(k^2b^{10})$ bit operations.
\end{abstract}

\section{Introduction}\label{sec:introduction}

Let $\PHI_k(N,M)$ be a random $k$-CNF formula on $N$ Boolean variables,
formed by choosing $M$ clauses uniformly without replacement from the
$2^k\binom Nk$ clauses on $k$ distinct variables. We refer to $r=M/N$ as the density. We say that an event occurs
\emph{with high probability} (w.h.p.) if its probability tends to 1 as
$N\to\infty$.

Recently, Carenini~\cite{Carenini} proved the satisfiability
threshold conjecture: for every $k\ge3$, there is a critical density,
$\rk$, below which a random formula is satisfiable w.h.p.\ and above which it is
unsatisfiable w.h.p. More than twenty years earlier, using the one-step
replica symmetry breaking (1RSB) cavity method, Mertens, M\'ezard and
Zecchina~\cite[\S7, Conjecture~1]{MMZ} predicted that $\rk$ equals the density
at which the complexity computed from the Survey Propagation (SP) equations
vanishes.

A rigorous approach to upper bounds on $\rk$ was developed by Franz and Leone~\cite{FL03} and Panchenko and Talagrand~\cite{PT04}. Their method, called \emph{interpolation}, gradually replaces clauses by auxiliary random terms, each depending on a single variable. The auxiliary terms are specified by probability distributions, called \emph{trial inputs}, which may be chosen freely within an admissible class. Given an admissible trial input the interpolation method  outputs a  lower bound on the expected minimum number of violated clauses. Suppose that at density $r$ this lower bound is at least $cN$ for some constant $c>0$ and all sufficiently large $N$. Concentration of the minimum number of violated clauses then implies that every assignment violates at least $cN/2$ clauses with exponentially high probability. In particular, $\rk\le r$.

In a major breakthrough, Ding, Sly and Sun (DSS)~\cite{DSS15} proved $\rk = \rstar$ for all $k\ge k_0$, where $k_0$ is an unspecified absolute constant.
Here $\rstar$ is the zero of the complexity obtained by iterating SP from
the point mass $\delta_{1/2}$, as described below.
Their upper bound uses SP to construct
a trial input for interpolation, while their matching lower bound uses an
SP-informed second-moment argument. Establishing the matching lower bound
constitutes the main technical contribution of their work~\cite[\S1.7]{DSS15}
and imposes the stronger constraints on $k_0$.
The upper-bound analysis also requires sufficiently large $k$ (at least
$3269$ for the estimates as written), and DSS note that they could not
extend their recursive SP analysis to all $k\ge3$.

Very recently, Vorobyev~\cite{Vorobyev} obtained $r_3\le4.268$ by choosing a trial input  tailored to $k=3$.

\subsection{Our contribution}

We prove $\rk \le \rstar$ for all $k \ge 3$, by extending the analysis of Ding, Sly and Sun~\cite{DSS15}.

The value $\rstar$ in the theorem of DSS comes from a particular
initialization of the SP recursion. (Throughout, as we will explain shortly, SP will be a recursion on probability
laws associated with the limiting random tree.) At each density $r$, DSS
start the SP recursion from the point mass $\delta_{1/2}$, so that the
initial survey value is always $1/2$. When the iteration converges,
its limiting law determines a quantity called ``complexity'' at that density. For all
sufficiently large $k$, DSS prove that the iteration converges throughout
an explicit density interval, that its limits satisfy the necessary
integrability conditions, and that the resulting complexity function has
a unique zero~\cite[Propositions~1.1--1.2 and Lemma~4.12]{DSS15}.

To extend this upper bound to every $k\ge3$, we first make explicit what the interpolation argument requires of a limiting SP law. Specifically, we state an explicit integrability assumption on the limiting law and verify that it justifies the expectations and limiting operations in the argument. We also handle the case of a nontrivial limiting law with zero complexity, without assuming convergence of SP at nearby densities. In that case, the argument proves unsatisfiability at every strictly larger density, which still gives an upper bound on the threshold at the original density. We collect this modest contribution as Theorem~\ref{prop:dss-if}. Our main contribution, Theorem~\ref{thm:uniform-sp}, is to prove that the SP iteration started from $\delta_{1/2}$ satisfies the requirements of this reduction for every $k\ge3$.

\paragraph{An upper bound on $\rk$ from SP properties (Theorem~\ref{prop:dss-if}, informal).}
Fix $k\ge3$ and a density $r$. If SP, started from any probability law on $[0,1)$, converges to a limit other than the one concentrated at zero, and this limit has finite expected logarithmic cost and nonpositive complexity, then $\rk\le r$. At each fixed density $s>r$, every assignment violates a positive linear number of clauses with exponentially high probability.

\paragraph{The required SP properties hold for every $k\ge3$ (Theorem~\ref{thm:uniform-sp}, informal).}
For every $k\ge3$ and every density $r$ in an explicit interval $H_k$,
the SP iteration started from $\delta_{1/2}$ converges to a nontrivial
limit with finite expected logarithmic cost. The complexity evaluated
at these limits has a unique zero $\rstar$ in that interval.
Applying Theorem~\ref{prop:dss-if} at this zero gives $\rk\le\rstar$
for every $k\ge3$.\bigskip

We also give certified enclosures of $\rstar$ for all $k\ge 3$ (Proposition~\ref{prop:enclosures}) and a method for efficient approximation to arbitrary precision (Theorem~\ref{thm:precision}). In Table~\ref{tab:trial-results} we give certified enclosures for $3 \le k \le 10$. The lower endpoints are lower bounds for $\rstar$, not the SAT threshold; the upper endpoints are  upper bounds for $\rstar$ and, thus, for $\rk$.

\begin{table}[!ht]\centering
\caption{Certified lower and upper bounds on the selected SP zero $\rstar$.}
\label{tab:trial-results}
\begin{tabular}{rrrrrr}
\toprule
$k$&Lower bound&Upper bound&$k$&Lower bound&Upper bound\\
\cmidrule(r){1-3}\cmidrule(l){4-6}
3&4.2667842934986&4.2667842940288&7&87.7852838240&87.7853002087\\
4&9.932385170&9.932389462&8&176.5425724542&176.5426073162\\
5&21.1172586283&21.1172625273&9&354.0102254068&354.0102263624\\
6&43.3698335799&43.3723251837&10&708.9152881881140&708.9152881896891\\
\bottomrule
\end{tabular}
\end{table}

After defining the relevant objects precisely, we state the reduction and our SP theorem
together, followed by their consequences. We then prove the reduction
to finite formulas. The remaining argument concerns SP: we give abstract
sufficient conditions for the required properties, verify them uniformly
for $k\ge10$, and then treat $k=3,\ldots,9$. The quantitative analysis
provides the certified enclosures and refinement algorithm.

\section{Survey Propagation and the upper-bound reduction}\label{sec:update}

To define the objects in the two main theorems, we first switch to an
auxiliary formula model. Let $\PHIaux_k(N,M)$ be formed by choosing all
$kM$ literal occurrences independently and uniformly from the $2N$ literals.
Variables may repeat within a clause, including with both signs, and clauses
may repeat. A standard conditioning argument transfers our upper bounds
from this model to the standard model $\PHI_k(N,M)$, preserving exponentially
high probability: condition first on distinct variables within every clause
and then on distinct clauses. The details are given in Section~\ref{sec:interp}.

We now pass from the auxiliary finite formula to its infinite-size local
limit. View the formula as a
signed bipartite factor graph, with variable vertices and clause vertices,
and explore it by breadth-first search from a uniformly chosen variable.
As $N\to\infty$ with $k$ fixed and $M/N\to r$, every fixed-depth
neighborhood converges in distribution to the corresponding neighborhood
of a rooted random tree \cite[\S4.1]{DSS15}. The SP recursion and its score
are defined on this limiting distribution of trees.

\subsection{The recursion on the limiting tree}\label{sec:sp-update}

Fix $k\ge3$ and $r>0$. At a variable of the limiting tree, distinguish
a literal $l$ and delete its parent clause, if present. Each remaining
incident clause has $k-1$ distinct child variables, and its edge to the
current variable carries one literal sign. Thus it is incident to either
$l$ or $\neg l$, never both.
A clause sends a \emph{warning} to its root literal when all its other
literals are forced to be false. A \emph{survey} is the probability
assigned by the recursion to the event that the remaining clauses force
a specified literal to be false.

Let $x$ be a probability law on $[0,1)$. For each of the two root literals,
independently sample a number of incident clauses with law $\Pois(kr/2)$,
where $\Pois(\lambda)$ denotes the Poisson law with mean $\lambda$.
For each clause take fresh independent surveys
$\surv_1,\ldots,\surv_{k-1}\sim x$, independently of the counts. Its warning
probability is $p=\prod_{j=1}^{k-1}\surv_j$. Conditional on these inputs,
the warnings are independent. Let $U$ be the product of $1-p$ over
clauses at $l$, and let $\overline U$ be the analogous product at
$\neg l$, with empty products equal to one.

For a given draw of the counts and input surveys, the probability that
warnings do not force both values is
$\mathfrak z=U+\overline U-U\overline U$. Conditional on this event, the
probabilities that $l$ is forced to be false, forced to be true, or left
free are respectively
\begin{equation}\label{eq:survey-triple}
 (\surv,\overline{\surv},F)=
 \frac{\big(U(1-\overline U),(1-U)\overline U,U\overline U\big)}
 {\mathfrak z}.
\end{equation}
Both $U$ and $\overline U$ are positive. The first coordinate therefore
lies in $[0,1)$, and its law over the original counts and input surveys
defines $T_r x$. The conditioning in~\eqref{eq:survey-triple} is performed
separately for each draw; it does not reweight the distribution of those
inputs. This construction is a distributional recursion on a tree. Its
use below as an interpolation input does not require identifying it with
message statistics of a finite formula.

Write $\delta_u$ for the point mass at $u$. The selected iteration is
\begin{equation}\label{eq:update}
 x^{(0)}=\delta_{1/2},\qquad
 x^{(n+1)}=T_r x^{(n)}\quad(n\ge0).
\end{equation}
At the first step, every clause has warning probability $2^{-(k-1)}$.
At each subsequent step, its inputs are independent draws from the law
just computed.

All logarithms are natural. Put $\Lfun(u)=-\ln(1-u)$ and call a scalar
survey law $x$ on $[0,1)$ \emph{admissible} if
$\E\Lfun(\surv)<\infty$ for $\surv\sim x$. For a variable update with
input $x$, let $\mathfrak z$ have the meaning above. The \emph{score} is
\begin{equation}\label{eq:dss-functional}
 \Gamma_r(x)=-\E\ln\mathfrak z
 -(k-1)r\,\E\Lfun\!\left(\prod_{j=1}^k\surv_j\right),
 \qquad \surv_1,\ldots,\surv_k\overset{\mathrm{iid}}\sim x.
\end{equation}
The two expectations are the variable cost and the clause cost.
Admissibility makes both costs finite and is preserved by $T_r$;
see Lemma~\ref{lem:basic-integrability} below. In
DSS~\cite[(1.9)--(1.11)]{DSS15}, $U,\overline U$ are denoted
by $\Pi^+,\Pi^-$ and $T_r$ by $\mathbf R$. At their selected limit, our
score is the negative of their $\Phi(r)$.

The density intervals used in our SP theorem are
\[
 H_3=[4.20,4.28],\qquad
 H_k=[2^k\ln 2-2,\,2^k\ln 2]\quad(k\ge4).
\]
For $k\ge4$, this is the interval of DSS~\cite[(1.7)]{DSS15}.

\subsection{The main results}\label{sec:sp-prediction}\label{sec:sp-obligations}

The first theorem is a conditional form of the upper-bound argument in
DSS~\cite[\S4.5, proof of Proposition~1.3]{DSS15}; we prove it in
Section~\ref{sec:interp}. Our SP theorem then supplies its hypotheses
at the selected zero.

\begin{theorem}[DSS upper-bound principle, conditional form]\label{prop:dss-if}
Fix $k\ge3$, $r>0$, and any probability law $x^{(0)}$ on $[0,1)$.
Let $x^{(n+1)}=T_r x^{(n)}$.
Suppose the iteration converges weakly to a probability law $\eta$
on $[0,1)$ such that
\[
 \eta\ne\delta_0,\qquad
 \int\Lfun(u)\,\eta(du)<\infty,\qquad
 \Gamma_r(\eta)\ge0.
\]
Then $\rk\le r$.  More precisely, for every fixed $s>r$ there are
$c,\gamma>0$ such that, for all sufficiently large $N$, with
probability at least $1-e^{-\gamma N}$ every assignment to
$\PHIaux_k(N,\lceil sN\rceil)$ violates at least $cN$ clauses.
The same holds for $\lfloor sN\rfloor$ uniform proper clauses,
sampled with or without replacement.
\end{theorem}

\begin{theorem}[The selected SP limit and its zero]\label{thm:uniform-sp}
For every integer $k\ge3$ and every $r\in H_k$, the iteration
\eqref{eq:update}, initialized at $\delta_{1/2}$, converges weakly to an
admissible fixed point $\eta_r\ne\delta_0$. The function $g_k(r)=\Gamma_r(\eta_r)$ is continuous and strictly increasing on $H_k$ and has a unique zero in
its interior. We denote this zero by $\rstar$ and call it the 1RSB
prediction.
\end{theorem}

The initialization in Theorem~\ref{prop:dss-if} is unrestricted;
convergence supplies a fixed point. Its score condition is nonpositive
complexity in the convention $\Phi=-\Gamma$.
At $r=\rstar$, Theorem~\ref{thm:uniform-sp} supplies a nontrivial
admissible limit with zero score, so Theorem~\ref{prop:dss-if} gives
$\rk\le\rstar$.

For sufficiently large $k$, DSS use the same map, initialization and
interval $H_k$. The iterates and their limits therefore agree, and the
selected zero is their prediction.

\begin{proposition}[Certified enclosures]\label{prop:enclosures}
For $3\le k\le10$, Table~\ref{tab:trial-results} gives strict lower
and upper bounds on $\rstar$. For every $k\ge10$, put
$t_k=2^k\ln 2-(1+\ln 2)/2$. Then
\begin{equation}\label{eq:asymptotic-window}
 t_k-\tfrac23 k^2 2^{-k}<\rstar
 <t_k-\tfrac1{10}k^2 2^{-k}.
\end{equation}
\end{proposition}

\begin{theorem}[Certified refinement]\label{thm:precision}
For all integers $k\ge3$ and $b\ge2$, a deterministic algorithm computes
a rational $\ell$ such that $\ell<\rstar<\ell+10^{-b}$ in $O(b^9\ln^2 b+k\ln^2 k)$ bit operations on a deterministic multitape
Turing machine, with an absolute implied constant. In particular, the
bit complexity is $O(k^2b^{10})$.
\end{theorem}

\begin{corollary}[The 1RSB upper bound]\label{cor:sat}
For every $k\ge3$, we have $\rk\le\rstar$. For each fixed
$r>\rstar$, there exist $c,\gamma>0$ such that, for all sufficiently
large $N$, with probability at least $1-e^{-\gamma N}$ every assignment
to $\PHI_k(N,\lceil rN\rceil)$ violates at least $cN$ clauses.
The same conclusion holds for the auxiliary model $\PHIaux_k(N,\lceil rN\rceil)$
and for uniform proper clauses sampled with replacement.
In particular, $r_3<4.2667842940288$.
\end{corollary}
\begin{proof}
At density $\rstar$, Theorem~\ref{thm:uniform-sp} supplies a nontrivial
admissible SP limit with zero score. Theorem~\ref{prop:dss-if}
therefore gives the threshold inequality and the finite-formula
conclusions. The numerical bound follows from the first row of
Table~\ref{tab:trial-results}.
\end{proof}

\subsection{The interpolation reduction}\label{sec:interp}\label{app:dss}

We now prove Theorem~\ref{prop:dss-if}. We use the interpolation inequality
and trial construction employed by DSS~\cite[\S4.5]{DSS15}. Their proof
applies the construction at strictly negative complexity. Here we verify
the limiting calculation for every $k\ge3$ under the stated logarithmic
integrability assumption. To handle a nontrivial SP fixed point of zero
complexity, we hold its law fixed and increase the density slightly; the
positive density derivative of its score then gives the strict inequality
required by interpolation. This uses only the fixed point at the original
density. The conclusion about a linear number of violated clauses with
exponentially high probability follows from the same bound and concentration.

We first justify the integrability assertions used in defining the SP
recursion and its score.

\begin{lemma}[Integrability]\label{lem:basic-integrability}
For every $k\ge3$ and $r>0$, an admissible input has an admissible
SP output and integrable variable and clause costs.
\end{lemma}
\begin{proof}\label{app:basic-integrability}
For either sign, the negative log no-warning probability is a
Poisson sum of terms $\Lfun(p)$, where $p$ is a product of $k-1$
input surveys. Since $p$ is at most its first factor, its expected
sum is at most $(kr/2)\E_x\Lfun(X)$. Now $\mathfrak z\ge U$ and
$1-\surv=\overline U/\mathfrak z\ge\overline U$, giving the same
bound for $\E[-\ln\mathfrak z]$ and $\E\Lfun(\surv)$.
The product bound also bounds the clause cost by $\E_x\Lfun(X)$.
\end{proof}

We now return to finite formulas, using the auxiliary model. Let $H_\Phi(\sigma)$ count
the clauses violated by an assignment $\sigma$, and put
\[
 Z_\beta(\Phi)=\sum_{\sigma\in\{0,1\}^N}e^{-\beta H_\Phi(\sigma)},
 \qquad
 f_{N,t}(\beta)=\frac1N\E\ln Z_\beta(\PHIaux_k(N,\boldsymbol M)),
 \qquad \boldsymbol M\sim\Pois(tN).
\]
The clause count is independent of the literal choices.  The following
auxiliary experiment defines the inputs to an interpolation bound on
$f_{N,t}$.

A bias $b\in(0,1)$ assigns probability $b$ to the value $1$.  For a
literal falsified by $e\in\{0,1\}$, write $q_0(b)=1-b$ and $q_1(b)=b$.
An \emph{interpolation input} is a probability law on probability
measures on $(0,1)$.  At each leaf, independently draw a measure from
this law and then a bias from that measure.  Take $D\sim\Pois(kt)$
clauses incident to one root, each with $k-1$ leaves and the root in its
last slot, and an independent clause indexed by $0$ with $k$ leaves.
Draw independent fair falsifying bits $e_{aj}$ in all slots, independently
of the measures and $D$.  For root value $s\in\{0,1\}$, set
\[
 C_a(s)=1-(1-e^{-\beta})\mathbf1\{s=e_{ak}\}
          \prod_{j=1}^{k-1}q_{e_{aj}}(b_{aj}),
 \qquad
 C_0=1-(1-e^{-\beta})\prod_{j=1}^kq_{e_{0j}}(b_{0j}).
\]
Here $\E_0$ averages over the biases conditional on the measures,
literal signs and $D$; the outer expectation averages over that
conditioning data.

\begin{proposition}[Interpolation bound]\label{prop:interpolation-bound}
For every $k\ge3$, $N\ge1$, $t,\beta>0$, every interpolation input and
$0<\theta<1$,
\begin{equation}\label{eq:dss-interpolation}
 \theta f_{N,t}(\beta)
 \le \E\ln\E_0\!\left[
       \left(\sum_{s=0}^1\prod_{a=1}^D C_a(s)\right)^\theta\right]
       -(k-1)t\,\E\ln\E_0[C_0^\theta].
\end{equation}
\end{proposition}
This is the bound of Panchenko and Talagrand~\cite[Theorem~3]{PT04},
following Franz and Leone~\cite{FL03}, in the form used by
DSS~\cite[\S4.5, Theorem~4.30]{DSS15}.  Its applicability to odd $k$ is
recorded in their footnote~18.  A field $y$ in that formulation
corresponds to the bias $b=e^y/(e^y+e^{-y})$.

\begin{lemma}[A positive-score interpolation input]\label{lem:survey-trial}
Let $(W_0,W_1,W_*)$ be a random probability vector whose law is invariant
under interchanging $W_0,W_1$, with admissible common marginal $x$.
Fix $k\ge3$ and $t>0$.  If $\Gamma_t(x)>0$, then there are $c,\gamma>0$
such that, for all sufficiently large $N$, with probability at least
$1-e^{-\gamma N}$ every assignment to
$\PHIaux_k(N,\lceil tN\rceil)$ violates at least $cN$ clauses.
The same holds for $\lfloor tN\rfloor$ uniform proper clauses,
sampled with or without replacement.
\end{lemma}
\begin{proof}
For $\beta>1$, put $q_\beta=(1+e^{2\beta})^{-1}$ and use independent
copies of the random measure
\[
 W_0\delta_{q_\beta}+W_1\delta_{1-q_\beta}+W_*\delta_{1/2}
\]
as the leaf measures in~\eqref{eq:dss-interpolation}.  Conditional on
the weights, sample labels $0,1,*$ independently.  A label $0$ or $1$
gives a bias favoring that value with probability $1-q_\beta$; the
label $*$ gives bias $1/2$.  Couple all $\beta$ using the same labels
and weights.  The conditional probability of a label forcing a
specified literal false has law $x$.

Conditional on the signs, weights and $D$, a clause warns against a
root value when its root literal is false and every leaf label forces
its literal false.  The probability that the root has an allowed value
is therefore $\mathfrak z$ from the SP update at density $t$ with
input $x$.  For the independent full clause, let $P$ be the product of
the $k$ forced-false masses; the probability that it has a leaf not
forced false is $1-P$.

For fixed labels, an allowed root value has
$\prod_a C_a(s)\ge2^{-D}$, whereas a forbidden one has an incident
factor at most $e^{-\beta}+(k-1)q_\beta$.  Similarly, $C_0\ge1/2$
unless all its labels force their literals false, in which case
$C_0\le e^{-\beta}+kq_\beta$.  These upper bounds are at most
$e^{-\beta/2}$ for all sufficiently large $\beta$, with $k$ fixed.
Taking $\theta=\beta^{-1/2}$ gives
\[
 \E_0\!\left[\left(\sum_s\prod_a C_a(s)\right)^\theta\right]
       \longrightarrow\mathfrak z,
 \qquad
 \E_0[C_0^\theta]\longrightarrow1-P.
\]
The powers are bounded by $2$ and $1$.  Their conditional expectations
have logarithms between $\ln\mathfrak z-D\ln2$ and $\ln2$, and
between $\ln(1-P)-\ln2$ and $0$, respectively.  The preceding
integrability bounds therefore justify passage to the outer
expectations.  The right side of~\eqref{eq:dss-interpolation} tends to
$-\Gamma_t(x)<0$.  Fixing sufficiently large $\beta$, we obtain
$\theta f_{N,t}(\beta)\le-\varepsilon$ for some $\varepsilon>0$ and
every $N$.

Let $V_N$ be the minimum violation count with the Poisson clause count.
An assignment attaining this minimum contributes $e^{-\beta V_N}$ to
$Z_\beta$, so $\E V_N\ge\varepsilon N/(\theta\beta)$.  Couple the
Poisson and deterministic clause counts by taking prefixes of one
sequence of independent clauses.  Their expected difference is at
most $\sqrt{tN}+1$, and adding or removing a clause changes the
minimum by at most one.  Thus the expected minimum remains linear at
either deterministic count.  Replacing a clause also changes the
minimum by at most one, so bounded differences gives the asserted
exponential bound in the independent-literal model.

For either deterministic clause count $m$, conditioning on distinct variables within
each clause gives independent uniform proper clauses.  This event
has probability
\[
 \left(\frac{(N)_k}{N^k}\right)^m
       \longrightarrow e^{-t\binom{k}{2}}>0.
\]
Conditioning further on distinct clauses gives sampling without
replacement.  Writing $L_N=2^k\binom Nk$, this second event has
conditional probability
$\prod_{j=0}^{m-1}(1-j/L_N)=1-O(N^{2-k})$.
Both conditionings preserve the exponential bound, after decreasing
$\gamma$ if necessary.
\end{proof}

We also need the density derivative with the input law held fixed.
For admissible $x$, write $x_h=\int u^h x(du)$ and
$d_h=\int u^h(T_t x)(du)$.  Adding a clause with warning probability
$p$ multiplies the consistency probability by $1-p\surv$.
The increment in the variable cost is $\Lfun(p\surv)$, bounded by
the integrable $\Lfun(p)$.  We may therefore differentiate the
Poisson series on each sign.  The two signs contribute
$k\E\Lfun(X_1\cdots X_{k-1}\surv)$, with independent
$X_i\sim x$ and $\surv\sim T_t x$; differentiating the explicit
clause term contributes $-(k-1)\E\Lfun(X_1\cdots X_k)$.
Expanding the logarithms yields
\begin{equation}\label{eq:density-derivative}
 \partial_t\Gamma_t(x)
 =\sum_{h\ge1}\frac{k d_hx_h^{k-1}-(k-1)x_h^k}{h}.
\end{equation}
The series is absolutely convergent by admissibility.

\begin{proof}[Proof of Theorem~\ref{prop:dss-if}]
Conditional on the Poisson counts, the SP update is a continuous
function of finitely many inputs in $[0,1)$: its denominator is
positive there.  Testing against bounded continuous functions and
then truncating the counts shows that $T_r$ is weakly continuous.
Hence the assumed limit satisfies $T_r\eta=\eta$.

One update at density $r$ with input $\eta$ gives the probability
triple~\eqref{eq:survey-triple}.  Its law is symmetric under exchanging
the two forced values, and both marginals equal $\eta$.
Hold this joint law fixed while varying the interpolation density.
By~\eqref{eq:density-derivative} and the fixed-point identity,
\[
 \left.\partial_t\Gamma_t(\eta)\right|_{t=r}
   =\sum_{h\ge1}\frac{\left(\int u^h\eta(du)\right)^k}{h}
   =\E\Lfun(Y_1\cdots Y_k)>0,
 \qquad Y_i\stackrel{\mathrm{iid}}\sim\eta.
\]
The expectation is finite by admissibility and positive because
$\eta\ne\delta_0$.  Since $\Gamma_r(\eta)\ge0$, for every $s>r$
there is a $t\in(r,s)$ with $\Gamma_t(\eta)>0$.
Apply Lemma~\ref{lem:survey-trial} at density $t$ to the fixed joint
law constructed at $r$.  Coupling each clause-sampling model by
adding clauses then gives the conclusion at density $s$.
Since this holds at every fixed $s>r$, $\rk\le r$.
\end{proof}

It remains to prove these SP results. Section~\ref{sec:abstract-sp}
gives abstract sufficient conditions and universal tools;
Sections~\ref{sec:uniform-regime} and~\ref{sec:finite-regime} verify them
for $k\ge10$ and $k=3,\ldots,9$. Section~\ref{sec:refine} gives the
refinement algorithm. The $k=10$ decimal calculation only sharpens
the uniform enclosure.

\section{The abstract SP argument}\label{sec:choice}\label{sec:abstract-sp}
The DSS criterion reduces the SAT upper bound to a nonzero admissible
fixed point of zero score. We state the SP conditions that suffice for
every $k\ge3$, then prove the analytic results used to verify them.

\subsection{Selection on a finite cover}\label{sec:orbit}
For a survey law $x$ write $x_h=\int v^h\,x(dv)$ and
$I(x)=\int(1-v)^{-1}\,x(dv)$. Positive majorants $u_h=O(h^{-1})$
define the moment metric
\begin{equation}\label{eq:moment-metric}
 D_u(x,y)^2=\sum_{h\ge1}\frac{u_h^{k-2}}h(x_h-y_h)^2.
\end{equation}
We write $D$ when the majorants are specified. Each contraction below
uses its declared majorants; different classes may use different metrics.

\begin{theorem}[Finite-cover criterion]\label{thm:finite-cover-sp}
\label{cor:abstract-zero}
Fix $k\ge3$. Let finitely many nondegenerate closed intervals
$J_i\subset H$ cover a nondegenerate compact interval
$H\subset(0,\infty)$. Suppose that for each $i$:
\begin{enumerate}[label=(C\arabic*)]
\item\label{cond:sp-classes}
$K_i$ is a nonempty convex class, weakly compact as laws on $[0,1]$,
with $I(x)\le B_i<\infty$ and $x_h\le u_{i,h}$ for every $x\in K_i$,
where $u_{i,h}>0$ and $u_{i,h}=O(h^{-1})$;
\item\label{cond:sp-dynamics}
for every $r\in J_i$, $T_rK_i\subseteq K_i$ and
$D_i(T_rx,T_ry)\le q_iD_i(x,y)$ for $x,y\in K_i$, where
$D_i=D_{u_i}$ and $0\le q_i<1$;
\item\label{cond:sp-entry}
for every $r\in J_i$, the initialized orbit $T_r^j\delta_{1/2}$
enters $K_i$ after finitely many steps.
\end{enumerate}
For every $r\in H$ the initialized orbit converges weakly to an
admissible law $\eta_r$, the unique fixed point and score maximizer
in each $K_i$ with $r\in J_i$.
The limits agree on overlaps. The selected score
$g_k(r)=\Gamma_r(\eta_r)$ is convex and continuously differentiable,
with derivative~\eqref{eq:selected-derivative} and one-sided endpoint
derivatives. The log-survey laws converge in Wasserstein distance of
order two.

If, in addition,
\textup{(P)} each $K_i$ has a common mean lower bound $\mu_i>0$,
then $\eta_r\ne\delta_0$ and $g_k$ is strictly increasing and strictly
convex. If also $r_-<r_+$ in $H$ satisfy $g_k(r_-)<0<g_k(r_+)$,
there is a unique zero in $H$, lying in $(r_-,r_+)$. At that zero
Theorem~\ref{prop:dss-if} applies. For $H=H_k$ the zero is $\rstar$.
\end{theorem}

\subsection{Convergence and score identities}

\begin{lemma}[Convergence in the moment metric]
\label{lem:abstract-convergence}
Fix $k\ge3$ and $r>0$. Let $K$ be a nonempty weakly compact class of
survey laws with $I(x)\le B_{\rm inv}<\infty$ and
$x_h\le u_h$, where $u_h>0$ and $u_h=O(h^{-1})$.
If $T_r(K)\subseteq K$ and $T_r$ contracts $D_u$ on $K$ by a factor
$q<1$, then $K$ is complete in $D_u$ and has a unique fixed law
$\eta$. For every $x\in K$, $T_r^jx$ converges weakly to $\eta$.
The same conclusion holds for the initialized orbit from $\delta_{1/2}$
if it enters $K$. In both cases the corresponding laws of the log survey
$\Lfun$ converge in Wasserstein distance of order two. The law $\eta$
is admissible; if $x_1\ge\mu>0$ throughout $K$, then
$\eta\ne\delta_0$.
\end{lemma}

\subsubsection{The two variations}\label{sec:score-changes}

\begin{proposition}[Variations of the score]\label{prop:score-variations}
Fix $k\ge3$ and a density $r>0$. Let $x,y$ be admissible survey laws,
put $d=T_rx$ and $\Delta_h=y_h-x_h$, and write $T=T_r$,
$\Gamma=\Gamma_r$. The fixed-input density derivative is
\eqref{eq:density-derivative}.
The right derivative $D\Gamma_x[y-x]$ of
$\Gamma((1-z)x+zy)$ at $z=0$ is
\begin{equation}\label{eq:variation}
 D\Gamma_x[y-x]
 =k(k-1)r\sum_{h\ge1}
   \frac{x_h^{k-2}(d_h-x_h)\Delta_h}{h}.
\end{equation}
The corresponding derivative formula applies at every interior point
of the segment. Both series are absolutely convergent. In particular,
every admissible fixed point is stationary in every admissible
mixture direction.
\end{proposition}

\begin{lemma}[Convexity in the warning intensity]
\label{lem:intensity-convexity}
Let $S(\mu)$ be the expected variable cost when the warning probabilities
on the two signs are independent Poisson point processes with the same
finite intensity measure $\mu$ on $[0,1)$, where
$\int\Lfun\,d\mu<\infty$. Then $S$ is convex on this class of
intensity measures. If $\mu((0,1))>0$, the function
$c\mapsto S(c\mu)$ is strictly convex for $c\ge0$.
Consequently, for an admissible survey law $x$, the function
$r\mapsto\Gamma_r(x)$ is convex, and it is strictly convex if
$x\ne\delta_0$.
\end{lemma}

\begin{lemma}[Score ascent]\label{lem:ascent}
Fix $k\ge3$ and $r>0$, and write $T=T_r$ and $\Gamma=\Gamma_r$.
For $s,t\in[0,1]$, set
$B_k(s,t)=s^k-t^k-kt^{k-1}(s-t)$.
For admissible survey laws $x,y$ and $d=Tx$,
\begin{equation}\label{eq:site-tangent}
 \Gamma(y)-\Gamma(x)
 \ge r\sum_{h\ge1}
       \frac{B_k(d_h,x_h)-B_k(d_h,y_h)}h.
\end{equation}
In particular $\Gamma(Tx)>\Gamma(x)$ unless $Tx=x$.
\end{lemma}
The Lean proofs of the ascent inequality and its strict consequence
assume finite inverse moments. The truncation argument in
Section~\ref{app:proof-score-ascent} extends them to the finite
logarithmic integrability assumed here.

\subsection{Regularity of the selected score}\label{sec:score}
\begin{proposition}[Regularity]\label{prop:regularity}
Fix $k\ge3$ and a nondegenerate compact density interval
$J\subset(0,\infty)$. Let $K$ be a nonempty weakly compact convex
class with the inverse and moment bounds of
Lemma~\ref{lem:abstract-convergence}. Suppose $T_rK\subseteq K$ and
$T_r$ contracts $D_u$ by a common factor $q<1$ for every $r\in J$.
Then its unique fixed law $\eta_r\in K$ depends weakly continuously
on $r$ and uniquely attains
\begin{equation}\label{eq:score-maximum}
 g(r)=\Gamma_r(\eta_r)=\max_{x\in K}\Gamma_r(x).
\end{equation}
The function $g$ is convex and continuously differentiable, with
one-sided endpoint derivatives and
\begin{equation}\label{eq:selected-derivative}
 g'(r)=\sum_{h\ge1}\frac{c_h(r)^k}{h},
 \qquad c_h(r)=\int v^h\,\eta_r(dv).
\end{equation}
If $x_1\ge\mu>0$ on $K$, then $g$ is strictly convex and
$g'(r)\ge\mu^k>0$.
\end{proposition}
The derivative is the expected clause cost
$\E\Lfun(\surv_1\cdots\surv_k)$ for independent surveys of law
$\eta_r$: independence gives $c_h(r)^k/h$ in its logarithmic expansion.

\subsection{Trial survey laws and endpoint signs}\label{sec:certificate}
A positive trial score in the same invariant class gives a positive selected score. To certify a
negative selected score, the trial score needs an upper correction.

\begin{proposition}[Trial score bounds]\label{prop:certificate}
Fix $t\in H_k$, put $T=T_t$, and let $\eta=\eta_t$ be an admissible
selected fixed law with moments $c_h(t)$. Let $K$ be a convex class of admissible laws
containing $x,\eta$, invariant under $T$, with moments at most
$u_h>0$, $u_h=O(h^{-1})$, and contraction factor $0\le q<1$ in
$D=D_u$. Write $\res_t(x)=D(Tx,x)$ and $C=\binom{k}{2}t$. Then
\begin{align}
 t B_k((Tx)_1,x_1)
 &\le g_k(t)-\Gamma_t(x)
   \le C\res_t(x)^2/(1-q),\label{eq:trial-budget}\\
 g_k(t)-\Gamma_t(x)
 &\le t\sum_{h\ge1}B_k(c_h(t),x_h)/h
   \le C D(\eta,x)^2.\label{eq:interval-comparison}
\end{align}
Suppose the selected branch on $[r_-,r_+]\subseteq H_k$ satisfies
\eqref{eq:selected-derivative}, with $c_1(r)\ge\mu>0$, and
$t\in[r_-,r_+]$. If $\gamma_-\le\Gamma_t(x)\le\gamma_+$ and
$g_k(t)-\Gamma_t(x)\le E$, then $g_k(r_-)<0<g_k(r_+)$ whenever
\begin{equation}\label{eq:certificate-margin}
 \gamma_++E<\mu^k(t-r_-),\qquad
 -\gamma_-<\mu^k(r_+-t).
\end{equation}
\end{proposition}

\subsection{Proofs of the abstract results}\label{app:survey-score-proofs}

\subsubsection{Convergence of the iteration}
\label{app:proof-survey-convergence}

\begin{proof}[Proof of Lemma~\ref{lem:abstract-convergence}]
The weights in~\eqref{eq:moment-metric} are positive and summable,
so $D$ convergence is equivalent to convergence of every moment,
and hence to weak convergence on $[0,1]$ by polynomial approximation.
The inverse cap gives $x([1-\varepsilon,1])\le B_{\rm inv}\varepsilon$.
For bounded continuous $f$ on $[0,1)$, extending it constantly from
$1-\varepsilon$ onwards changes its integral by at most
$2B_{\rm inv}\varepsilon\|f\|_\infty$. Thus weak convergence on
$[0,1]$ also gives weak convergence on $[0,1)$ within $K$.
Compactness makes $K$ complete in $D$. Banach's theorem gives the
unique fixed law and convergence of every orbit in $K$, including
the initialized orbit after entry. The inverse cap implies
admissibility; the positive mean bound, when assumed, passes to the
limit and excludes $\delta_0$.

For the initialized orbit, $\Lfun(\surv)\le-\ln\overline U$ is
bounded by a Poisson sum of costs $\Lfun(p)$, each at most
$\Lfun(\surv_1)$ for an incoming survey in that clause. Second
moments are therefore finite by induction from $\delta_{1/2}$.
On $K$, $\Prb(\Lfun(Y)>z)\le B_{\rm inv}e^{-z}$ makes the squared
logarithms uniformly integrable. Continuity of $\Lfun$ gives weak
convergence of their laws, and uniform integrability gives convergence
of their second moments, hence convergence in Wasserstein distance
of order two, or mean square convergence under a suitable coupling.
\end{proof}

\subsubsection{The two variations of the score}
\label{app:proof-score-variations}

\begin{proof}[Proof of Proposition~\ref{prop:score-variations}]
The density derivative was proved in Section~\ref{sec:interp}.
Write $\Pi_j$ for the product of $j$ independent surveys with law $x$,
independent of the current output $\surv$.
For variation in the input law, the derivative is
\[
 k(k-1)r\int\!
 \big[\E\Lfun(v\Pi_{k-2}\surv)
        -\E\Lfun(v\Pi_{k-1})\big]\,(y-x)(dv).
\]
Conditional on the counts, differentiation replaces one input law at
a time by the signed measure $y-x$. In the variable term, first subtract
the cost with the affected clause removed; this is permitted because
$y-x$ has mass zero. Each remaining increment is bounded by
$\Lfun(v)$, which is integrable under $x+y$, uniformly along the
segment. The number of replacements is $k-1$ times the incident clause
count and has finite expectation, so the conditional derivatives can
be averaged. For a Poisson count $D$, the identity
$j\Prb(D=j)=\E[D]\Prb(D=j-1)$ restores the original law of the
other clauses. Expanding the logarithms proves~\eqref{eq:variation}.
Its absolute series is bounded by a constant multiple of
$\sum_h(x_h+y_h)/h$. The same argument works at every interior point
of the mixture segment.
\end{proof}

\subsubsection{Convexity in the warning intensity}
\label{app:proof-intensity-convexity}

\begin{proof}[Proof of Lemma~\ref{lem:intensity-convexity}]
Split each sign process into two independent halves of intensity
$\mu/2$. In one half let $\surv,\overline{\surv}$ be the
conditional probabilities of forcing the test literal false and true,
and put $m_{ij}(\mu)=\E[\overline{\surv}^{\,i}\surv^j]$.
Sign symmetry gives $m_{ij}=m_{ji}$. Use primes for an independent
half. Given consistency of both halves, their union conflicts with
probability $\overline{\surv}\surv'+\surv\overline{\surv}'$,
so joining them adds
$-\ln(1-\overline{\surv}\surv'-\surv\overline{\surv}')$ to their
log costs. For a second intensity $\nu$, expand this logarithm for
halves of intensities $(\mu/2,\mu/2)$, $(\nu/2,\nu/2)$ and
$(\mu/2,\nu/2)$. Independence and symmetry give the moment factors
$m_{ij}(\mu)^2$, $m_{ij}(\nu)^2$ and
$m_{ij}(\mu)m_{ji}(\nu)=m_{ij}(\mu)m_{ij}(\nu)$.
The separate half-costs cancel, leaving
\begin{equation}\label{eq:split-root}
 S(\mu)+S(\nu)-2S((\mu+\nu)/2)
 =\sum_{i+j\ge1}\frac1{i+j}\binom{i+j}{i}
       [m_{ij}(\mu)-m_{ij}(\nu)]^2.
\end{equation}
Each nonnegative logarithmic expansion has finite sum, justifying
subtraction and proving midpoint convexity.

For $\mu_z=(1-z)\mu+z\nu$ and $0\le s<t\le1$, couple the
processes by sharing intensity $(1-t)\mu+s\nu$ on each sign.
The remaining intensities are $(t-s)\mu$ and $(t-s)\nu$.
An added warning $p$ changes the variable cost by at most $\Lfun(p)$,
so $|S(\mu_t)-S(\mu_s)|\le2(t-s)\int\Lfun\,d(\mu+\nu)$.
This continuity turns midpoint convexity into convexity.

On a ray with nonzero warnings, adding warnings decreases
$F=(U^{-1}+\overline U^{-1}-1)^{-1}$ strictly with positive
probability. For a half of intensity $c\mu/2$, its expectation is
$1-2m_{10}(c\mu)$ and strictly decreases with $c$. Thus distinct
$c,c'$ give a positive $(1,0)$ term in~\eqref{eq:split-root}, proving
strict convexity along the ray. For input $x$ the sign intensity
$(kr/2)\Law(\Pi_{k-1})$, with $\Pi_{k-1}$ the product of $k-1$
independent $x$-surveys, varies on such a ray when $x\ne\delta_0$.
The clause term is affine in $r$, proving the final assertion.
\end{proof}

\subsubsection{Score ascent}\label{app:proof-score-ascent}

\begin{proof}[Proof of Lemma~\ref{lem:ascent}]
For input law $x$, each sign has warning intensity
$\mu_x=(kr/2)\Law(\Pi_{k-1})$. Apply
Lemma~\ref{lem:intensity-convexity} along the straight segment from
$\mu_x$ to $\mu_y$. The single-warning calculation gives its
directional derivative at $\mu_x$ as
\[
 kr\sum_{h\ge1}d_h(y_h^{k-1}-x_h^{k-1})/h.
\]
Convexity bounds $S(\mu_y)-S(\mu_x)$ below by this derivative.
Subtracting the change of the clause term gives
\eqref{eq:site-tangent}, first for inputs bounded away from one. The
algebra is
\[
 B_k(d_h,x_h)-B_k(d_h,y_h)
 =(k-1)(x_h^k-y_h^k)+kd_h(y_h^{k-1}-x_h^{k-1}).
\]

For general admissible $x,y$, truncate every input survey $Y$ to
$\min(Y,1-\varepsilon)$, using the same samples and counts. Each
update uses finitely many values, all strictly below one; hence its
truncated and original outputs agree for all sufficiently small
$\varepsilon$ on each realization. Every output moment converges.
The truncated variable and clause costs are bounded by the original
$-\ln U$ and $\Lfun(\surv_1)$, respectively, for an original
clause input $\surv_1$. These dominators are integrable.

In the difference of the two $B_k$ terms, the common $d_h^k$ cancels.
Truncation decreases the input moments, and every output moment is
at most one. The absolute value of the resulting $h$th summand before
the factor $r$ is therefore at most
$(2k-1)(x_h+y_h)/h$. This is summable because
$\sum_h x_h/h=\int\Lfun(v)x(dv)$ and likewise for $y$.
Dominated convergence proves~\eqref{eq:site-tangent} in full generality.
Taking $y=d$ leaves the nonnegative terms $B_k(d_h,x_h)$.
Since $B_k(s,t)>0$ for $s\ne t$, equality is possible only if all
moments agree, which by moment determinacy means $Tx=x$.
\end{proof}

\begin{lemma}[Positive score]\label{lem:dss-field}
Suppose that $k\ge3$, $t>0$, $y$ is admissible and $d=T_t y$.
If either $\Gamma_t(y)>0$ or $\Gamma_t(d)>0$, the conclusion of
Lemma~\ref{lem:survey-trial} holds at density $t$.
\end{lemma}
\begin{proof}
Run one step of the recursion at $t$ with input $y$ for the positive
literal of a variable. The three output probabilities
$(\surv,\overline{\surv},F)$ give the probabilities of false, true and
joker. Interchanging the signs of the incident clauses exchanges the
first two entries and preserves their joint law. Both entries have law
$d$, which is admissible. By Lemma~\ref{lem:ascent},
$\Gamma_t(d)\ge\Gamma_t(y)$, so either hypothesis gives
$\Gamma_t(d)>0$. Apply Lemma~\ref{lem:survey-trial} to this triple.
\end{proof}

\subsubsection{Regularity on one class and compatibility across a cover}
\label{app:proof-regularity}

\begin{proof}[Proof of Proposition~\ref{prop:regularity}]
Lemma~\ref{lem:abstract-convergence} supplies the unique fixed law in
$K$ at each density. Since
$I(x)=\sum_{i\ge0}x_i\ge(h+1)x_h$, every $x\in K$ has
$x_h\le B_{\rm inv}/(h+1)$. The logarithmic tail bound in the proof
of that lemma gives uniformly bounded second moments of $\Lfun$.

The update is jointly weakly continuous in the density and input law.
To see this, couple counts at nearby densities by adding independent
Poisson counts at the larger density; the probability of any added
clause tends to zero. For fixed counts, the output survey is a bounded
continuous function of finitely many input surveys in $[0,1)$.
Weak convergence of their product laws therefore gives weak convergence
of the output; averaging over the counts proves the assertion.
The variable cost is bounded by $-\ln U$, a Poisson sum with uniformly
bounded second moment on the compact density interval. The clause cost
is bounded by $\Lfun(\surv_1)$ for any one input. Thus both costs
are uniformly integrable, and $\Gamma_r(x)$ is jointly continuous.
Because $T_rx\in K$, the density derivative series is uniformly
dominated by
\[
 (2k-1)\sum_{h\ge1}\frac{[B_{\rm inv}/(h+1)]^k}{h},
\]
and hence its sum is jointly continuous as well.

Compactness gives a score maximizer in $K$. Invariance and strict ascent
force every maximizer to be fixed; uniqueness of the fixed law gives
\eqref{eq:score-maximum}. If $r^{(n)}\to r$ in $J$ and a subsequence of
$\eta_{r^{(n)}}$ converges weakly to $y$, then for every $x\in K$ the
inequality $\Gamma_{r^{(n)}}(\eta_{r^{(n)}})\ge\Gamma_{r^{(n)}}(x)$ passes to
the limit. Thus $y$ maximizes $\Gamma_r$ and equals $\eta_r$.
Compactness now proves weak convergence of the whole sequence, so the
fixed law depends weakly continuously on $r$.

For $s<t$ in $J$, put $\Delta\Gamma(x)=\Gamma_t(x)-\Gamma_s(x)$.
Using $\eta_s$ as a trial survey law at $t$ and $\eta_t$ as a trial
survey law at $s$ gives
\[
 \Delta\Gamma(\eta_s)\le g(t)-g(s)\le\Delta\Gamma(\eta_t).
\]
Divide by $t-s$. Each outer term is a difference quotient with its law
held fixed, and the mean value theorem makes it a partial density
derivative at an intermediate density. Joint continuity of that derivative
and weak continuity of the fixed law give the same limit on both sides
as $s,t\to r$, namely $\partial_r\Gamma_r(\eta_r)$.
The limit is one-sided at the endpoints. Substitution of $d_h=x_h$
in~\eqref{eq:density-derivative} proves~\eqref{eq:selected-derivative}
and continuity of the derivative. If $x_1\ge\mu>0$ on $K$, its
first term is at least $\mu^k$.

For each fixed $x\in K$, Lemma~\ref{lem:intensity-convexity} makes
$\Gamma_r(x)$ convex in $r$, and strictly convex if $x\ne\delta_0$.
For $s<v<t$, apply this convexity to the particular law $x=\eta_v$
and bound its scores at $s,t$ by $g(s),g(t)$. This proves convexity
of $g$, and proves strict convexity under the positive-mean assumption.
\end{proof}

\begin{proof}[Proof of Theorem~\ref{thm:finite-cover-sp}]
On each $J_i$, Lemma~\ref{lem:abstract-convergence} and the entry
condition give the initialized limit in $K_i$, including log-$W_2$
convergence. At overlaps both are limits of the same sequence
$T_r^j\delta_{1/2}$, so they agree even when the metrics differ.
Proposition~\ref{prop:regularity} gives maximization and regularity on
each interval. Refine $H$ at the finitely many panel endpoints.
The laws are continuous on each piece and agree at its joins; the
one-sided score derivatives agree there by
\eqref{eq:selected-derivative}. Thus $g_k$ is $C^1$ on $H$.
Its derivative is nondecreasing on each piece and across the joins,
so $g_k$ is convex. Under \textup{(P)} the derivative is positive and
strictly increasing on each piece, hence strictly increasing on $H$.
Finally the opposite endpoint signs give a zero by continuity, and
strict increase makes it unique. Its law is nonzero and admissible,
so Theorem~\ref{prop:dss-if} applies.
\end{proof}

\subsubsection{Bounds from a trial survey law}
\label{app:proof-trial-score}

\begin{proof}[Proof of Proposition~\ref{prop:certificate}]
For $x_\theta=(1-\theta)x+\theta Tx$, equation~\eqref{eq:variation}
and Cauchy--Schwarz give
\[
 \left|\frac{d}{d\theta}\Gamma_t(x_\theta)\right|
 \le2C\res_t(x_\theta)\res_t(x),\qquad
 \res_t(x_\theta)\le D(Tx_\theta,Tx)+D(Tx,x_\theta)
 \le(q\theta+1-\theta)\res_t(x).
\]
The four laws $x,Tx,x_\theta,Tx_\theta$ lie in $K$. Their moments are
at most $u_h$, so the $h$th derivative summand is bounded absolutely
by $2C u_h^k/h$, uniformly in $\theta$. Integration gives
\[
 \Gamma_t(Tx)-\Gamma_t(x)
 \le2C\res_t(x)^2\int_0^1(1-\theta+q\theta)\,d\theta
 =C(1+q)\res_t(x)^2.
\]
Contraction gives $\res_t(T^jx)\le q^j\res_t(x)$ and
$D(T^jx,\eta)\le q^jD(x,\eta)$. The latter proves weak convergence
to the given fixed law directly, without a completeness assumption.

The scores converge as well. For $h\ge2$ and any $Y$ with law in $K$,
\[
 \Prb(Y>1-1/h)\le\frac{u_h}{(1-1/h)^h}=O(h^{-1}).
\]
Taking $h=\lfloor e^z\rfloor$ gives
$\Prb(\Lfun(Y)>z)=O(e^{-z})$, uniformly on $K$.
Thus the second logarithmic moments are uniformly bounded; the
Poisson-sum bound for the variable cost and the one-input bound for
the clause cost give uniform integrability as in the proof of
Proposition~\ref{prop:regularity}. The score is consequently continuous
along these iterates.

Summing their gains proves the upper bound in~\eqref{eq:trial-budget},
because $(1+q)\sum_{j\ge0}q^{2j}=1/(1-q)$.
For the lower bound, monotone ascent together with convergence of the
scores gives $\Gamma_t(Tx)\le\Gamma_t(\eta)=g_k(t)$.
Take $y=Tx$ in~\eqref{eq:site-tangent} and retain its first nonnegative
term. This proves the other side of~\eqref{eq:trial-budget}.

Apply~\eqref{eq:site-tangent} at $\eta$, with comparison law $x$.
Since $T\eta=\eta$, it gives
$g_k(t)-\Gamma_t(x)\le t\sum_hB_k(c_h(t),x_h)/h$.
For $a,b\in[0,u]$, Taylor's theorem gives
$B_k(a,b)\le\binom{k}{2}u^{k-2}(a-b)^2$, proving
\eqref{eq:interval-comparison}.
Finally,
$\gamma_-\le g_k(t)\le\gamma_++E$ and
$g_k'(r)\ge\mu^k$ on $[r_-,r_+]$ by the stipulated derivative
formula. Integrating the slope on the two sides of $t$ gives the strict
endpoint signs under~\eqref{eq:certificate-margin}.
\end{proof}

\subsection{Common coordinates and evaluation identities}\label{sec:common-coordinates}
For a survey $X$ drawn from $x$ write $c_h=x_h=\E X^h$ for its moment of order $h$. Set $c_0=1$, $m=k-1$ and
$\lambda=kr/2$. The warning probability of a fresh clause has the law of $p_\circ=\prod_{j=1}^m \surv_j$ for
independent $\surv_j\sim x$. Write $S=-\ln U$ and $\bar S=-\ln\overline U$ for the two log fields. They are independent
Poisson sums of copies of $Y=\Lfun(p_\circ)$, each with mean count $\lambda$.
It is convenient to keep track of the field difference $Z$ and of ${\mathcal V}$,
the exponential of minus half their sum. In these coordinates we use
\begin{samepage}
\[
 \surv=\frac{e^{\bar S}-1}{e^{S}+e^{\bar S}-1},\quad
 F=\frac1{e^{S}+e^{\bar S}-1},\quad Z=\bar S-S,\quad
 {\mathcal V}=e^{-(S+\bar S)/2},\quad \sigma_Z=(1+e^{-Z})^{-1}.
\]
Here $F=U\overline U / \mathfrak z$ is the probability that the variable is a joker, conditional on the absence of a conflict.
The joint law is symmetric in $S,\bar S$. Exchanging them sends $\surv=\sigma_Z-F(1-\sigma_Z)$ to $\overline{\surv}=1-F-\surv$,
so $\E \surv=(1-\E F)/2$.
\end{samepage}

Each invariant class is defined by bounds on the moments $c_h$ and on
the inverse moment $I(x)=\E_x(1-X)^{-1}\le B_{\rm inv}$.
These constraints are convex and weakly closed on laws on $[0,1]$ by moment
continuity and lower semicontinuity of $(1-t)^{-1}$, set to $+\infty$ at
$1$. Hence the classes are weakly compact and exclude mass at $1$.
For $h\ge2$ the inverse moment bound gives
\begin{equation}\label{eq:inverse-caps-common}
 c_h\le\mathcal U(h,B_{\rm inv}):=(B_{\rm inv}-1)(h-1)^{h-1}/h^h\le (B_{\rm inv}-1)/(2h).
\end{equation}
Indeed $\E[X/(1-X)]\le B_{\rm inv}-1$ and $t^{h-1}(1-t)$ is maximized at $(h-1)/h$. Since $I(x)=\sum_{i\ge0}c_i$, summing the AM--GM inequalities $t^i+t^{2h-i}\ge2t^h$ for $j\le i<h$ and adding $t^h$ gives
\begin{equation}\label{eq:head-tail}
 (2h-2j+1)c_h\le I(x)-\sum_{i<j}c_i,\qquad h\ge j .
\end{equation}

\subsubsection{Universal evaluation identities}
\begin{lemma}[Three two-field identities]\label{lem:fourier-identities}
For every \(k\ge3\), every \(r>0\), and every admissible survey law,
the site cost, output moments, and warning-response operator have the
following representations.
For a function $C(u,v)$ of the two
log fields, inserting a warning $p$ on either sign gives the operator
\[
 \Delta_p C=C(u+a,v)+C(u,v+a)-2C(u,v),\qquad a=-\ln(1-p).
\]
For the variable cost $C=-\ln\mathfrak z$, sign symmetry gives
$\E\Delta_p C=2\sum_{\ell\ge1}(T_r x)_\ell p^\ell/\ell$. Let $B_{\ell h}$ be the
coefficient of $p^h$ in $\E\Delta_p(\surv^\ell)$. Thus the row indexed by $\ell$
records how the $\ell$th output moment responds to an extra warning as
its probability varies. Put $\mathsf K_{\ell h}=\sqrt{h/\ell}\,B_{\ell h}$ for $\ell,h\ge1$. Use
$\widehat f(\xi)=\int_{\mathbb R} f(y)e^{-i\xi y}\,dy$ and the branch of $\ln(1-p)$
vanishing at zero, with radial boundary limits.
\par\Needspace{11\baselineskip}
For mode $n\ge0$, frequency $\xi\ge0$
and $\alpha=n/2+i\xi$, put
\begin{gather*}
 a_h(\alpha)=(-1)^h\binom{\alpha}{h},\qquad f_{n,h}(\xi)=\sqrt h\,\Re a_h(\alpha),\qquad
 g_n(\xi)=\frac{2|\Gamma_{\mathrm E}(\alpha)|^2}{\pi n!},\\
 \Psi_n(\xi)=\E[{\mathcal V}^ne^{i\xi Z}]=\exp\{2\lambda[\Re\E(1-p_\circ)^{\alpha}-1]\}.
\end{gather*}
Here $\Gamma_{\mathrm E}$ is Euler's gamma function. For $n=0$ the kernel is defined for $\xi>0$. At $\xi=0$ we use the limits of the full integrands given in Section~\ref{sec:fourier}. The following identities give the
variable cost, the updated moments and their response to an inserted warning.
Write $f_n=(f_{n,h})_{h\ge1}$ and use $\preceq$ for the positive
semidefinite order on every finite block.
\begin{align}
 \E C&=\lambda\E Y-\ln 2+
   \tfrac12\sum_{n\ge0}\int_0^\infty g_n(\xi)[\Psi_n(\xi)-\mathbf1_{n=0}]\,d\xi,\label{eq:fourier-site}\\
 (T_r x)_\ell&=\tfrac12+\tfrac{\ell}{2}\sum_{n\ge0}\int_0^\infty \Psi_n g_n\,\Re a_\ell(\alpha)\,d\xi,\label{eq:fourier-moments}\\
 \mathsf K&=\sum_{n\ge0}\int_0^\infty \Psi_n g_n f_nf_n^{\mathsf T}\,d\xi,\qquad
 0\preceq\mathsf K\preceq I.\label{eq:three-gram}
\end{align}

\end{lemma}
\begin{proof}
See Section~\ref{sec:fourier} for the transform calculation and its
truncation and discretization errors.
\end{proof}

\section{The uniform range \texorpdfstring{$k\ge10$}{k at least 10}}\label{sec:uniform-regime}\label{app:selected-sp}
The estimates below close the argument for $k\ge10$.
The class and contraction lemmas also apply
for $k\ge5$ and $k\ge6$, respectively; their remaining scalar bases are
verified in Section~\ref{sec:finite-scalar-bases}.
\begin{proposition}\label{prop:uniform-obligations}
For $k\ge10$, the SP iteration from $\delta_{1/2}$ converges throughout
$H_k$ to an admissible nonzero fixed law $\eta_r$. Its score $g_k$ is
continuously differentiable and strictly convex, with derivative
\eqref{eq:selected-derivative} and a unique zero in the interval
\eqref{eq:asymptotic-window}.
\end{proposition}
The Fourier/Gram identities used below are derived in
Section~\ref{sec:fourier}.

\subsection{A common invariant class}\label{app:common}
Since $(1-\surv)^{-1}=1+(e^{\bar S}-1)e^{-S}$ and the two fields are
independent, their Poisson transforms give
\begin{equation}\label{eq:inverse-identities}
 I(T_r x)=1+\exp\big(\lambda\textstyle\sum_{h\ge2}c_h^m\big)-e^{-\lambda c_1^m}.
\end{equation}

For every $k\ge5$ put $m=k-1$, $\rho=2^{-k}$, $\mu=(1-\rho)/2$, $\epsilon=1.21k\rho$ and $B_{\rm inv}=11/4$. Impose
\begin{equation}\label{eq:invariant-class}
 \E X\ge\mu,\qquad \E X^h\le b_h:=\epsilon/2+(1-\epsilon)2^{-h}\ (h\ge1),
 \qquad \E(1-X)^{-1}\le B_{\rm inv},
\end{equation}
Let $K_k$ be the class of laws satisfying these inequalities. It is convex
and weakly compact by the inverse cap, and it contains the starting law
$\delta_{1/2}$. The inverse cap gives the moment bounds
\begin{equation}\label{eq:invariant-caps}
 u_1=1/2,\qquad u_h=\min\{b_h,\mathcal U(h,B_{\rm inv})\}\quad(h\ge2).
\end{equation}
For convergence and endpoint comparisons we use the moment majorants
\begin{equation}\label{eq:checked-common-metric}
 \widehat u_h=\min\{b_h,1/h\}\quad(h\ge1),\qquad
 \widehat D_k(x,y)^2=\sum_{h\ge1}\frac{\widehat u_h^{\,k-2}}h
                         (x_h-y_h)^2 .
\end{equation}
Both sequences majorize the moments on $K_k$: for $h\ge2$ the inverse
cap gives $u_h\le7/(8h)<1/h$, while $u_1=\widehat u_1=1/2$.
We use $u_h$ for scalar estimates and $\widehat D_k$ for convergence and
endpoint comparisons. The contraction proof applies to each metric separately.

\begin{lemma}[An invariant class for $k\ge5$]\label{lem:invariance}
For $k\ge5$ and $r\in H_k$, $T_rK_k\subseteq K_k$.
With $\lambda_-=k(2^k\ln2-2)/2$ and $\lambda_+=k2^k\ln2/2$,
we have $z=e^{-\lambda_-\mu^m}<2\rho$.
\end{lemma}
\begin{proof}
Put $C=(B_{\rm inv}-1)/2=7/8$,
$B=\sum_{h\ge2}2\mathfrak h_{h-1}u_h^m/h$ and
$J=\sum_{h\ge2}u_h^m$. Expanding the square of
$Y=-\ln(1-p_\circ)$ gives $\E Y^2\le B$.

\emph{Mean.} Put $z_1=z\exp(-\lambda_-\mu^{2m}/4)$. The Poisson transform gives
\[
\E {\mathcal V}\le z_1,\qquad \E {\mathcal V}^2\le z^2,\qquad \E {\mathcal V}^3\le z^3\exp(\lambda_+u_2^m).
\]
Here $\sqrt{1-u}\le1-u/2-u^2/8$ and
$(1-u)^{3/2}\le1-3u/2+u^2/2$ at $u=p_\circ$. Sign exchange and $F\le {\mathcal V}/(2-{\mathcal V})\le {\mathcal V}/2+({\mathcal V}^2+{\mathcal V}^3)/4$ now give
\begin{equation}\label{eq:mean-invariance}
 \E \surv\ge\tfrac12-z_1/4-\big[z^2+z^3\exp(\lambda_+u_2^m)\big]/8 .
\end{equation}

\emph{Moment profile.} Since $\surv\le \sigma_Z=(1+\tanh(Z/2))/2$, expanding $\sigma_Z^h$ and bounding every positive even power of $\tanh$ by its square gives $\E \surv^h\le (1-d)2^{-h}+d/2$ with $d=\E\tanh^2(Z/2)$, since odd powers have zero expectation.
It suffices to show $d\le\epsilon$. Write $v=\lambda_+ B$. The characteristic function of $Z$ is $\exp[2\lambda(\E\cos(\xi Y)-1)]\ge e^{-v\xi^2}$. For $f(z)=\operatorname{sech}^2(z/2)$, the transform
$\widehat f(\xi)=4\pi \xi/\sinh(\pi \xi)$ is nonnegative and integrable.
Writing $\varphi_Z$ for the characteristic function above, Fourier inversion gives
\[
 1-d=\frac1{2\pi}\int_{\mathbb R}\widehat f(\xi)\varphi_Z(\xi)\,d\xi
 \ge\frac1{2\pi}\int_{\mathbb R}\widehat f(\xi)e^{-v\xi^2}\,d\xi
 =\E\operatorname{sech}^2 X_{\rm G},
 \qquad X_{\rm G}\sim\mathcal N(0,v/2),
\]
since $e^{-v\xi^2}$ is the characteristic function of $2X_{\rm G}$. By Cauchy--Schwarz, for $a=4/5$,
\begin{equation}\label{eq:gaussian-comparison}
 1-d\ge\frac{(\E e^{-aX_{\rm G}^2})^2}{\E[e^{-2aX_{\rm G}^2}\cosh^2 X_{\rm G}]}
 =\frac{2\sqrt{1+2av}}{(1+av)[1+\exp(v/(1+2av))]} .
\end{equation}
Alternatively $\tanh^2 y\le y^2$ gives $d\le\E Z^2/4\le v/2$, since
$\E Z^2=2\lambda\E Y^2$ by independence of the two Poisson sums.

\emph{Inverse moment.} Since $c_1\le1/2$, \eqref{eq:inverse-identities} gives $I(T_r x)\le1+e^{\lambda_+J}-e^{-\lambda_+2^{-m}}$.
The cases $k=5,6,7$ are verified in Section~\ref{sec:finite-scalar-bases}.

\emph{Uniform verification for $k\ge8$.} We have $B\le J$ since $2\mathfrak h_{h-1}\le h$. We have $4u_2\le1+\epsilon$, $u_h<0.14$ for $h\ge3$, and $m\epsilon<0.265$. Thus
\[
 4^m J\le e^{0.265}+16C^2(0.56)^{m-2}\sum_{h\ge3}h^{-2}<1.7 .
\]
It follows that $v/2\le1.7(\ln2)k\rho<\epsilon$ and $1+e^{\lambda_+J}<11/4$. Finally $\lambda_-\mu^m=k(\ln2-2\rho)(1-\rho)^m\ge k\ln2-k\rho(2+m\ln2)$ by Bernoulli's inequality. The loss $k2^{-k}(2+(k-1)\ln2)$ is below $\ln(1.7)$ at $k=8$ and decreases thereafter. The successive ratios of $k2^{-k}$ and $k(k-1)2^{-k}$ are $(k+1)/(2k)$ and $(k+1)/(2(k-1))$, respectively. Hence $z<1.7\rho$, and \eqref{eq:mean-invariance} gives $\E \surv>(1-\rho)/2$ using $z_1\le z$ and $e^{\lambda_+u_2^m}<2$.
\end{proof}

\subsection{The common contraction estimate}\label{app:moment-contraction}
\begin{lemma}[Analytic contraction for $k\ge6$]\label{lem:contract}\label{lem:common-contraction}
For every $k\ge6$ and every $r\in H_k$, the SP update on $K_k$ has
Lipschitz factor at most $q_k$ of \eqref{eq:common-factor}, both in the
metric formed from \eqref{eq:invariant-caps} and in
$\widehat D_k$ of \eqref{eq:checked-common-metric}. Here
$q_6<0.797$, $q_{k+1}<q_k$ for $k\ge6$, and $q_k<1/3$ for $k\ge8$.
\end{lemma}
\begin{proof}
The insertion and block reduction hold for $k\ge3$ on any
inverse-moment-bounded class with positive moment caps $a_h=O(1/h)$.
For this lemma take $a_h=u_h$ or $\widehat u_h$, write $D_a$ for its
metric, and put $w_h=a_h^{m-1}$, $W_0=\operatorname{diag}(w_h)$.
The calculation applies separately to each choice. For a segment $x_\theta=(1-\theta)x+\theta y$ in the class, put $c_h=(x_\theta)_h$ and $\Delta c_h=y_h-x_h$.
The warning moment $c_h^m$ has derivative $m c_h^{m-1}\Delta c_h$.
Inserting a warning at rate $\lambda$ on each sign gives
\begin{equation}\label{eq:moment-insertion}
 \frac{d}{d\theta}(T_r x_\theta)_\ell=m\lambda\sum_{h\ge1} B_{\ell h}c_h^{m-1}\Delta c_h .
\end{equation}
Since $\overline{\surv}=1-F-\surv$, inserting a warning of probability $p$ gives the two possible output surveys
$(\surv+Fp)/(1-\overline{\surv}p)$ and $\surv(1-p)/(1-\surv p)$. For $|p|\le1$,
$|\surv+Fp|\le \surv+F=1-\overline{\surv}\le|1-\overline{\surv}p|$ and
$|1-\surv p|^2-\surv^2|1-p|^2=(1-\surv)(1+\surv-2\surv\Re p)\ge0$.
Cauchy's bound on their $\ell$th powers gives $|B_{\ell h}|\le2$ for $h\ge1$,
and the moment tails justify differentiation. In the coordinates
$\sqrt{w_h/h}\,\Delta c_h$, the derivative is
$m\lambda W_0^{1/2}\mathsf K W_0^{1/2}$ times a diagonal operator of
norm at most one, since $0\le c_h^{m-1}/w_h\le1$. Integrating along the
segment converts the uniform derivative bound into a Lipschitz bound.

For $\beta\ge m\lambda$, $\tau>0$, put
$\Lambda=\operatorname{diag}((\beta w_h-\tau)_+)$. Since $w_h\to0$,
$\Lambda$ has finite rank. Conjugating $\beta W_0\preceq\tau I+\Lambda$
by $\mathsf K^{1/2}$, using $0\preceq\mathsf K\preceq I$ and
$\|B^*B\|=\|BB^*\|$, gives
\begin{equation}\label{eq:moment-block}
 \|T_r'\|_{D_a}\le\tau+\|\Lambda^{1/2}\mathsf K \Lambda^{1/2}\|.
\end{equation}
Thus the finite matrix with $d_h=(\beta w_h-\tau)_+>0$ bounds the
operator in the full moment metric.
\paragraph{A three-coordinate bound for \texorpdfstring{$k\ge6$}{k ≥ 6}.}
Take $\beta=0.694k(k-1)/(2\rho)\ge m\lambda_+$. The operator
$\mathsf K^{1/2}W_0\mathsf K^{1/2}$ increases in positive semidefinite
order with $W_0$ and has the same norm as $W_0^{1/2}\mathsf K W_0^{1/2}$.
Since both cap sequences are bounded by $b_h$, we may use
$w_h=b_h^{k-2}$ in the operator majorant in \eqref{eq:moment-block},
retaining the chosen metric $D_a$ on its left. This proves the bound
in each metric directly. Taking $\tau=\beta w_4$ leaves only
coordinates $1,2,3$. Put $z=2\rho$. Lemma~\ref{lem:invariance}
and its proof give $\Psi_1\le\E {\mathcal V}\le z$ and $\Psi_n\le\E {\mathcal V}^2\le z^2$
for $n\ge2$. Restrict $f_n$ and $\mathsf K$ to the first three coordinates.
Bounding the sum of the undamped Gram matrices for $n\ge2$ by $I$ in
\eqref{eq:three-gram} gives
\[
 \mathsf K_{[3]}\preceq \mathsf G_0+z\mathsf G_1+z^2 I,
 \qquad \mathsf G_j=\int_0^\infty g_j(\xi)f_j(\xi)f_j(\xi)^{\mathsf T}\,d\xi .
\]
Put $a=(0,1,\sqrt{3/2})$, $b=(0,1,\sqrt{3/8})$ and $e_1=(1,0,0)$. Directly from the first three binomial coefficients, $f_0=-\xi^2 a/\sqrt2$ and $f_1=-e_1/2-(\xi^2+1/4)b/\sqrt2$. Hence $\mathsf G_0=aa^{\mathsf T}/8$ and $\mathsf G_1\preceq(e_1e_1^{\mathsf T}+bb^{\mathsf T})/2$. Indeed $\int \xi^4g_0=1/4$, and $g_1$ has mass one and moments $1/4,5/16$ at orders two and four, by the cosine transform above. Use $(u+v)^2\le2u^2+2v^2$ for $\mathsf G_1$.

On coordinates $2,3$ bound the rank-one terms by their trace and use $d_h\le \beta w_h$ in \eqref{eq:moment-block}. The diagonal term $z^2 I$ on those coordinates contributes at most $\beta z^2w_2$. This gives the rational contraction factor
\begin{equation}\label{eq:common-factor}
 \|T_r'\|_{D_a}\le q_k:=\beta \left[w_4+\max\left\{(z/2+z^2)w_1,\,
 (1/8+z/2)(w_2+3w_3/2)+z^2w_2\right\}\right].
\end{equation}
Exact rational evaluation gives $q_8<1/3$; the $k=6$ base is checked in
Section~\ref{sec:finite-scalar-bases}. When $k$ increases by one, $\beta $ is multiplied by $2(k+1)/(k-1)$,
while $w_1=4\rho$ and $z=2\rho$ both halve. Since $\epsilon$ decreases,
$b_h(k+1)\le b_h(k)\le b_2(k)<0.28$ for $h\ge2$, so
$w_h(k+1)/w_h(k)\le b_h(k)<0.28$. Each summand inside the brackets
therefore has successive ratio at most $0.28$ (for $zw_1$ it is $1/4$).
The same bound holds for their sums and maximum, giving
$q_{k+1}/q_k\le0.56(k+1)/(k-1)<1$ for $k\ge6$.

\end{proof}

\subsection{The two endpoint signs}\label{app:high}
The class and contraction lemmas give the selected branch on $H_k$
and its derivative \eqref{eq:selected-derivative}. We now establish the
opposite signs needed for its zero.

\begin{lemma}[Uniform initial enclosure]\label{lem:uniform-endpoints}
For $k\ge10$, put $\rho=2^{-k}$, ${\ell_2}=\ln2$,
$d=(1+{\ell_2})/2$, $t={\ell_2}/\rho-d$, and $x=\delta_{1/2}$.
The two densities $r_-=t-\tfrac23 k^2\rho$ and
$r_+=t-\tfrac1{10}k^2\rho$ lie in $H_k$ and satisfy
$g_k(r_-)<0<g_k(r_+)$. Thus
\begin{equation*}
 t-\tfrac23 k^2\rho<\rstar<t-\tfrac1{10}k^2\rho.
\end{equation*}
\end{lemma}
\begin{proof}
For $x=\delta_{1/2}$, every clause sends a warning of probability
$(1/2)^{k-1}=2\rho$.
Write $v=kt$, $\omega=-\ln(1-2\rho)/2$, so
$U=e^{-2\omega D_0}$ and $\overline U=e^{-2\omega D_1}$, where $D_0,D_1$
are the independent Poisson counts of mean $v/2$ at $l,\neg l$.
Put $X=\omega(D_1-D_0)$, $Y=e^{-\omega(D_0+D_1)}$ and $w=Y/(2\cosh X)$.
Since $\mathfrak z=2Y\cosh X(1-w)$ and $-\omega+\ln\cosh\omega=\ln(1-\rho)$,
\begin{equation}\label{eq:point-score}
 -\Gamma_t(x)={\ell_2}+t\ln(1-\rho)+\mathcal D+\E\ln(1-w),\qquad \mathcal D=\E\ln\cosh X-v\ln\cosh\omega .
\end{equation}
For the mixed moments, put $m_j=\E Y^j$ and $\nu_j=v\omega^2e^{-j\omega}$. Independence of
$D_0,D_1$ gives
$\E[Y^je^{sX}]=\exp\{v(e^{-j\omega}\cosh(s\omega)-1)\}$.
Its second and fourth derivatives at $s=0$ give
$\E[Y^jX^2]=m_j\nu_j$ and $\E[Y^jX^4]=m_j(\nu_j\omega^2+3\nu_j^2)$.

\emph{The trial score.} Expanding \eqref{eq:point-score} gives
\begin{equation}\label{eq:signed-point-score}
 \left|\Gamma_t(x)-\rho^2\left[-\tfrac18+\tfrac{{\ell_2}}{12}
             +\tfrac{1-{\ell_2}}{4}k+\tfrac{{\ell_2}^2}{4}k^2\right]\right|\le2k^3\rho^3 .
\end{equation}
Indeed, $\rho\le\omega\le \rho+\rho^2/(1-2\rho)$, $\nu_0<k\rho$ and $|\nu_0-k{\ell_2}\rho|\le3k\rho^2$. Put $\zeta=1-\sqrt{1-2\rho}=\rho+\zeta^2/2$. Then $\rho\le\zeta\le \rho+\rho^2$ and $|v\zeta-k{\ell_2}+k\rho/2|\le2k\rho^2$. Substituting in $m_1=e^{-v\zeta}$, $m_2=e^{-2v\rho}$ and $m_3=e^{-v(\zeta+2\rho-2\rho\zeta)}$ gives
\[
 |m_1-\rho-k\rho^2/2|\le k^2\rho^3,\quad |m_2-\rho^2|\le2k\rho^3,\quad
 m_3\le2\rho^3,\quad |m_1\nu_1-k{\ell_2}\rho^2|\le2k^2\rho^3 .
\]
These use $k\rho\le10/1024$ and $|e^u-1-u|\le u^2e^{|u|}/2$.

The real Taylor bound $0\le\ln\cosh z-z^2/2+z^4/12\le z^6/45$ follows by integrating $z^2-2z^4/3\le\tanh^2 z\le z^2$ twice in $(\ln\cosh z)''=1-\tanh^2 z$. Since $\E X^6=\nu_0\omega^4+15\nu_0^2\omega^2+15\nu_0^3$, it gives $|\mathcal D+(k{\ell_2}\rho)^2/4|\le k^3\rho^3/2$. Likewise $|\operatorname{sech}z-1+z^2/2|\le2z^4$ by the bound $\sup_{z\in\mathbb R}|\operatorname{sech}^{(4)}z|\le45$. Expanding $-\ln(1-w)$ through $w^2/2$, with $w\le1/2$, gives
\[
 \left|-\E\ln(1-w)-\left(m_1/2-m_1\nu_1/4+m_2/8\right)\right|
 \le\E[YX^4]+\E[Y^2X^2]/8+m_3/12\le4k^2\rho^3 .
\]
Finally $|{\ell_2}+t\ln(1-\rho)-\rho/2-(d/2-{\ell_2}/3)\rho^2|\le \rho^3$. Combining these bounds proves \eqref{eq:signed-point-score}.

\emph{The limiting score.} Put $y=T_t x$ and
$\Delta=1/2-y_1=\E[w/(1-w)]/2$, using sign symmetry.
Since $w\le Y/2\le1/2$, the bounds on $m_1,m_2$ give
$0\le\Delta\le(m_1+m_2)/4<\rho/3$. Using $\operatorname{sech}X=1+O(X^2)$ in $w$ gives
$\Delta=\rho/4+O(k\rho^2)$. Expand $[(1+\tanh X)/2]^j$, which bounds $\surv^j$ from above since $\surv=(U-U\overline U)/\mathfrak z\le U/(U+\overline U)$.
Symmetry cancels odd powers and every nonconstant even power of $\tanh X$ is at most
$\tanh^2 X\le X^2$, so $y_j-2^{-j}\le\E X^2/2<k\rho/2$.
Jensen's inequality gives $2^{-j}-y_j\le j2^{1-j}\Delta$.
Set $C=\binom k2t$ and $S_{\rm tail}=\sum_{j\ge2}\widehat u_j^{k-2}/j$. Equation~\eqref{eq:trial-budget} gives
\begin{equation}\label{eq:gain-sandwich}
 C y_1^{k-2}\Delta^2\le g_k(t)-\Gamma_t(x)
 \le\frac C{1-q}\left[4\rho\Delta^2+(k\rho)^2 S_{\rm tail}/4\right].
\end{equation}
Here $S_{\rm tail}\le 0.26^{k-2}/2+0.55\cdot 0.13^{k-3}$ since $\widehat u_2<0.26$ and $\widehat u_j\le\min(0.13,1.1/j)$ thereafter. Using $q=1/3$, the upper gain divided by $k^2\rho^2$ is at most $(3{\ell_2}/4)[4/9+k^2 S_{\rm tail}/(4\rho)]<0.4$. This bound decreases with $k$.

Equation~\eqref{eq:signed-point-score} places the trial score between $0.105k^2\rho^2$ and $0.15k^2\rho^2$. Write $c_j(r)=\int v^j\,\eta_r(dv)$. Throughout $H_k$, $0.98\rho<\mu^k\le g_k'\le \rho+0.875(0.26)^{k-1}<1.01\rho$ by the derivative formula, $\mu^k\ge \rho(1-k\rho)$, $c_1(r)\le1/2$ and $\sum_{j\ge1}c_j(r)\le7/4$. By \eqref{eq:gain-sandwich},
$0.105k^2\rho^2<g_k(t)<(0.15+0.4)k^2\rho^2$. Both endpoints in
\eqref{eq:asymptotic-window} lie in $H_k$. Integrating the derivative bounds from the endpoints to $t$ gives a negative
score at the left endpoint since $0.15+0.4<(2/3)\cdot 0.98$, and a positive
score at the right endpoint since $0.105>0.1\cdot1.01$.

\end{proof}

\begin{proof}[Proof of Proposition~\ref{prop:uniform-obligations}]
The compact convex class $K_k$ contains $\delta_{1/2}$ and has positive
mean and a finite inverse cap. Lemmas~\ref{lem:invariance} and
\ref{lem:contract} give invariance and contraction in $\widehat D_k$;
Lemma~\ref{lem:uniform-endpoints} gives opposite endpoint signs.
Theorem~\ref{thm:finite-cover-sp} proves the claim, and
Theorem~\ref{prop:dss-if} gives the threshold upper bound.
\end{proof}

\subsection{The large-\texorpdfstring{$k$}{k} correction}
The same calculation identifies the next asymptotic term.
Using $q=q_k$ in \eqref{eq:gain-sandwich}, the two sides agree to $o(\rho^2)$ as $k\to\infty$. Indeed $y_1^{k-2}=4\rho(1+O(k\rho))$, $k^4S_{\rm tail}/\rho\to0$, and the three-coordinate bound in Section~\ref{app:moment-contraction} gives $q_k=O(k^2(0.56)^k)$. Thus $g_k(t)-\Gamma_t(x)$ is ${\ell_2}k(k-1)\rho^2/8+o(\rho^2)$. The derivative bounds likewise give $g_k'/\rho=1+o(k^{-2})$ uniformly on $H_k$. The mean value identity $\rstar=t-g_k(t)/g_k'(\widetilde r)$, with $\widetilde r$ between $\rstar$ and $t$, gives
\begin{equation}\label{eq:cavity-polynomial}
 \rstar=t+2^{-k}\left[\tfrac18-\tfrac {\ell_2}{12}
        +\tfrac{3{\ell_2}-2}{8}k-\tfrac{{\ell_2}+2{\ell_2}^2}{8}k^2\right]+o(2^{-k}).
\end{equation}
Mertens, M\'ezard and Zecchina~\cite[Eq.~(A.22)]{MMZ} obtained the same
expansion as $k\to\infty$.

\section{The remaining finite range}\label{sec:finite}\label{sec:finite-regime}
\label{app:numerical}\label{app:trials}

For $6\le k\le9$ the common class and three-coordinate contraction
leave only scalar base checks and endpoint evaluations. At $k=5$ two
spectral tests supply contraction. At $k=4,3$ we also construct invariant
classes and prove entry of the initialized orbit. We first state the finite
endpoint bounds, then treat these cases in order.

\begin{lemma}[Finite endpoint certificates]\label{lem:table-check}
For every $3\le k\le10$, the trial survey law specified in this section
belongs to the invariant class at both certificate densities $r_-,r_+$
in Table~\ref{tab:endpoint-margins}. There are certified bounds
\[
 u_-\ge\Gamma_{r_-}(x),\qquad
 \varepsilon\ge\mathcal D_k(T_{r_-}x,x)^2,\qquad
 l_+\le\Gamma_{r_+}(x)
\]
such that
\[
 u_-+\binom{k}{2}\frac{r_-\varepsilon}{1-q}<0,
 \qquad l_+>0.
\]
Here $q=0.98$ for $k=3$, $q=0.99$ for $k=4,5$, and $q=q_k$
from~\eqref{eq:common-factor} for $k\ge6$.
For $k=3,4,5$, $\mathcal D_k$ uses the checked moment caps of
Proposition~\ref{prop:finite-contraction}; for $k\ge6$ it is
$\widehat D_k$ of~\eqref{eq:checked-common-metric}.
\end{lemma}

By the trial score bounds, the first inequality gives $g_k(r_-)<0$.
The second gives $g_k(r_+)\ge\Gamma_{r_+}(x)>0$. The lower bounds below
include the stated error allowances and are rounded down to two
significant figures.

\begin{table}[!ht]\centering
\small
\setlength{\tabcolsep}{4pt}
\caption{Original certificate densities and certified score margins.
The fourth and fifth columns are positive lower bounds; the last gives
an upper bound for $g_k'$ throughout $[r_-,r_+]$.}
\label{tab:endpoint-margins}
\begin{tabular}{rrrrrr}
\toprule
$k$&$r_-$&$r_+$&$-g_k(r_-)$&$\Gamma_{r_+}(x)$&$M_k$\\
\midrule
3&4.2667842934&4.26678429403&$2.0\cdot10^{-11}$&$2.5\cdot10^{-13}$&0.20281673048863\\
4&9.93238&9.93239&$4.9\cdot10^{-7}$&$5.1\cdot10^{-8}$&0.09476443926736\\
5&21.1172328074&21.1172637876&$8.4\cdot10^{-7}$&$4.1\cdot10^{-8}$&0.03253168788971\\
6&43.3691402552&43.3724008191&$1.1\cdot10^{-5}$&$1.2\cdot10^{-6}$&0.01586558125752\\
7&87.7850548417&87.7853574542&$1.8\cdot10^{-6}$&$4.5\cdot10^{-7}$&0.00786086926362\\
8&176.5424396848&176.5426583813&$5.2\cdot10^{-7}$&$2.0\cdot10^{-7}$&0.00391656426708\\
9&354.0102059738&354.0102330105&$3.8\cdot10^{-8}$&$1.3\cdot10^{-8}$&0.00195543446586\\
10&708.9152881482&708.9152882644&$3.9\cdot10^{-11}$&$7.3\cdot10^{-11}$&0.00097709921496\\

\bottomrule
\end{tabular}
\end{table}

\begin{lemma}[Sharpening the zero enclosures]\label{lem:enclosure-tightening}
The selected fixed laws in the panel or common invariant classes satisfy
$g_k'(r)<M_k$ on the original certificate interval. Their zero
therefore satisfies the refined bounds in Table~\ref{tab:trial-results}.
\end{lemma}
\begin{proof}
The endpoint densities $r_-,r_+$ in Table~\ref{tab:endpoint-margins}
are the inputs to the preceding certificates. Write $A_-,A_+$ for the
tabulated positive lower bounds on $-g_k(r_-)$ and $g_k(r_+)$.
If $g_k'(r)\le M_k$ throughout this interval, integration gives
\begin{equation}\label{eq:refined-zero-bracket}
 r_-+A_-/M_k\le\rstar\le r_+-A_+/M_k.
\end{equation}
We obtain the slope bounds in the same table from the existing moment
bounds, as follows. Take $(\mu,B_{\rm inv})=(2/5,14/3)$ for $k=3$,
$(9/20,5)$ for $k=4$, and $((1-2^{-k})/2,11/4)$ otherwise.
Write $c_h=c_h(r)$. At the selected fixed point, symmetry gives $c_1\le1/2$ and Jensen's
inequality gives $c_j\ge\mu^j$. Cauchy--Schwarz gives
$c_h\le c_{h-1}-c_{h-1}^2/(I(\eta_r)-\sum_{j=0}^{h-2}c_j)$.
Starting with $v_1=1/2$, for $2\le h\le256$ compute
\[
 D_h=B_{\rm inv}-\sum_{j=0}^{h-2}\mu^j,\qquad
 v_h=\min\{u_h,\ v_{h-1}-v_{h-1}^2/D_h\},
\]
rounding each $v_h$ upwards to a multiple of $10^{-18}$.
Here $u_h=\mathcal U(h,B_{\rm inv})$ for $k=3,4$ and the caps
\eqref{eq:invariant-caps} otherwise. The substitution is monotone
because $D_h\ge2v_{h-1}$: for $k=3,4$ we have $D_h\ge3$; otherwise
$D_2=7/4$ and, for $h\ge3$, $D_h\ge3/4$ and $v_{h-1}<0.30$.
Thus $c_h\le v_h$, uniformly over the original bracket. By
\eqref{eq:selected-derivative} and \eqref{eq:inverse-caps-common},
\[
 g_k'(r)\le\sum_{h=1}^{256}\frac{v_h^k}{h}
 +\frac{[(B_{\rm inv}-1)/2]^k}{k\,256^k}<M_k,
\]
where the displayed $M_k$ is rounded upwards to $10^{-14}$.
Rounding the endpoints in~\eqref{eq:refined-zero-bracket} strictly
outwards gives Table~\ref{tab:trial-results}.

\end{proof}

\subsection{Two-atom witnesses for \texorpdfstring{$6\le k\le10$}{6 through 10}}
\label{sec:finite-two-atom}

The common-class estimates apply directly for $8\le k\le10$; the scalar
checks below extend them to $k=6,7$. Each endpoint witness has two atoms
of equal mass $1/2$:

\begin{table}[htp]\centering
\caption{Two-atom trial survey laws and the densities used by the
independent comparison checks. All displayed numbers are exact decimals.}
\label{tab:two-atom-laws}
\begin{tabular}{rrrr}
\toprule
$k$&$t$&First atom&Second atom\\
\midrule
6&43.3724008191&0.3649&0.6254\\
7&87.7852&0.3934&0.6021\\
8&176.54255&0.4216&0.5762\\
9&354.01022&0.44246&0.55651\\
10&708.915288205&0.45773422&0.54176100\\
\bottomrule
\end{tabular}
\end{table}

We evaluate the same survey law at both $r_-$ and $r_+$. The test density
$t$ is used in~\eqref{eq:certificate-margin}; it equals $r_+$ for $k=6$
and lies inside the bracket for $k=7,8,9,10$. The head caps $b_h$ and
tail cap $1/h$ bound the caps of $\widehat D_k$ in the residual calculation.

\subsubsection{The finite scalar bases}\label{sec:finite-scalar-bases}

We record here the scalar calculations needed for $k=5,6,7$, so that
the $k=5$ case below can use the same common class.
Retain $B,J,\lambda_-,\lambda_+$ from the proof of
Lemma~\ref{lem:invariance}, with $\lambda_-,\lambda_+$ as defined there
from the endpoints of $H_k$.
For each of these three values, sum $B$ and $J$ through $H=128$
using the caps in~\eqref{eq:invariant-caps}.
With $C=(B_{\rm inv}-1)/2$, the bound $u_h\le C/h$ and the integral
test give tails at most
\[
 2(C/H)^m\left(\frac{\mathfrak h_H}{m}+\frac1{m-1}\right),
 \qquad
 \frac{C^mH^{1-m}}{m-1},
\]
respectively. Substitution gives a right side greater than $\mu$ in
\eqref{eq:mean-invariance}, greater than $1-\epsilon$ in
\eqref{eq:gaussian-comparison}, and
\[
 1+e^{\lambda_+J}-e^{-\lambda_+2^{-m}}<B_{\rm inv}.
\]
It also gives $e^{-\lambda_-\mu^m}<2\rho$. These finite scalar
inequalities complete common-class invariance at $k=5,6,7$.
Evaluating~\eqref{eq:common-factor} at $k=6$ gives $q_6<0.797$;
the successive-ratio estimate in its proof gives the required
contraction also at $k=7$. The same scalar certificate contains the
numerical base estimates at $k=8$, including $q_8<1/3$, used in the
uniform argument.

For the trial survey laws at $5\le k\le10$, rational arithmetic verifies
the inverse cap, mean and first eight moment bounds in
\eqref{eq:invariant-class}. All subsequent moment bounds follow from
\[
 x_h\le x_8\le\epsilon/2\le b_h\qquad(h>8).
\]
The endpoint scores and residuals are evaluated in
Section~\ref{sec:finite-endpoint-evaluation}.

\subsection{The exceptional contraction estimate at \texorpdfstring{$k=5$}{k=5}}
\label{sec:finite-five}

For the common class at $k=5$, two $14$-coordinate spectral tests at
the ends of $[20.18,22.19]\supset H_5$ give contraction by $0.99$
throughout the interval. The same construction handles $k=4,3$.

\begin{proposition}[Contraction in the three remaining cases]
\label{prop:finite-contraction}
For the common class at $k=5$ and every $r\in H_5$, the SP update
has Lipschitz factor at most $0.99$.
For $k=4,3$, every panel class $K_J$ of Lemma~\ref{lem:entry}, and
every $r\in J\cap H_k$, the factors are at most $0.99$ and $0.98$,
respectively. Each assertion uses the moment metric formed from the
checked caps described in Section~\ref{sec:finite-spectral}.
At $k=3$ a single metric is used on each group of ten consecutive panels;
at $k=4$ each panel has its own metric.
\end{proposition}

The formal checks use these declared caps for $k=3,4,5$ and
$\widehat D_k$ for $k\ge6$.

The five-atom trial survey law is specified by the following integer
pairs. Both the atom and mass in each pair are multiplied by $10^{-8}$.

\begin{table}[htp]\centering
\caption{The five-atom trial survey law for $k=5$.}
\label{tab:five-atom-law}
\begin{tabular}{rr}
\toprule
Atom&Mass\\
\midrule
12677294&5028632\\
29069273&26069141\\
49382210&39795528\\
69797437&24648431\\
86497642&4458268\\
\bottomrule
\end{tabular}
\end{table}

\subsection{Panel classes for \texorpdfstring{$k=4$ and $k=3$}{k=4 and k=3}}
\label{sec:finite-small}\label{app:class-entry}

For $k=4,3$ use one class for each closed density subinterval, called
a panel. Moment enclosures prove invariance and entry of the initialized
orbit; the spectral tests then give contraction.

\begin{lemma}[Entry and invariance for $k=3,4$]\label{lem:entry}
Partition $[9.09,11.091]\supset H_4$ into sixteen equal closed intervals,
and $[4.20,4.28]=H_3$ into forty equal closed intervals.
Each panel $J$ has a weakly compact convex class $K_J$ such that
$T_rK_J\subseteq K_J$ for every $r\in J\cap H_k$.
For every $r\in J\cap H_k$, the iteration from $\delta_{1/2}$ satisfies
$x^{(j)}\in K_J$ for all $j\ge31$ when $k=4$ and all $j\ge49$
when $k=3$. The $k=4$ classes satisfy $x_1>9/20$ and $I(x)\le5$;
the $k=3$ classes satisfy $x_1>2/5$ and $I(x)\le14/3$.
\end{lemma}

Each class has $64$ rational moment intervals, an inverse cap and a
strictly feasible atomic witness. Section~\ref{sec:panel-certificate}
gives the construction. Intersecting the panels with $H_k$ yields a
nondegenerate finite cover, trimming the two end panels at $k=4$.
At shared endpoints the same initialized orbit determines both local limits.

The endpoint trial survey laws have eight atoms at $k=4$ and ten at
$k=3$. In each pair, multiply both entries by $10^{-5}$ for $k=4$
and by $2^{-26}$ for $k=3$.

\begin{table}[htp]\centering
\caption{Trial survey laws for the two panel regimes.}
\label{tab:small-atom-laws}
\setlength{\tabcolsep}{4pt}
\begin{tabular}{rrrrrrrr}
\toprule
\multicolumn{4}{c}{$k=4$}&\multicolumn{4}{c}{$k=3$}\\
\cmidrule(r){1-4}\cmidrule(l){5-8}
Atom&Mass&Atom&Mass&Atom&Mass&Atom&Mass\\
\midrule
4055&2384&57667&21642&2214661&9691394&60561008&3783639\\
12480&9999&73499&15859&12133005&16350431&64218844&1630229\\
25169&18210&86477&7729&27496576&15768874&65992790&563145\\
40862&22599&95336&1578&42597683&11846422&66762001&140552\\
&&&&53755640&7315682&67042208&18496\\
\bottomrule
\end{tabular}
\end{table}

Rational arithmetic checks all $64$ moment bounds and the inverse cap
in panel $7$ for $k=4$, and panel $34$, namely $[4.266,4.268]$, for
$k=3$. Both endpoints of each zero bracket lie in the indicated panel.
The independent comparison at $k=4$ uses $t=r_+$.
At $k=3$ the left certificate requires moments through order $8192$;
the right certificate proves positivity at the sharper exact density
$4.26678429403$.

\subsection{Soundness of the panel certificates}\label{sec:panel-certificate}
\begin{proof}[Proof of Lemma~\ref{lem:entry}]
For each panel $J$, the 64 pairs of rational moment
bounds form its \emph{box}. Its class $K=K_J$ consists of the laws whose
moments lie in that box and that satisfy $I(x)\le B_{\rm inv}$, where $B_{\rm inv}=14/3$ for $k=3$ and $B_{\rm inv}=5$ for $k=4$. The stored orders are $1$ through $24$, then steps of
$2,4,8,16$ up to $48,96,192,256$, respectively.
The moment bounds are integer multiples of $10^{-9}$. A strictly feasible
rational atomic law witnesses nonemptiness of each class.

We use the bounds in two ways. For invariance the input can be any law
in the box satisfying the inverse-moment bound. For entry we follow the
particular orbit starting at $\delta_{1/2}$ until its enclosure lies
inside the box and satisfies that bound.

Let $s_h=c_h$ or $c_h^m$, the moments of an input survey or of a clause warning probability, respectively.
Between successive stored orders $i,j$, let $P_h$ interpolate $s_i,s_j$
linearly. For a
stored order $i_0\le i/2$ and $i\le h\le j$,
\begin{equation}\label{eq:moment-secant}
 0\le P_h-s_h\le\frac{(h-i)(j-h)}{2(i-i_0)^2}s_{i_0} .
\end{equation}
For fixed $t=e^{-y}$, the function $f(u)=t^u$ is convex, so its chord
on $[i,j]$ lies above it. At $h$ the gap is at most
$\tfrac12(h-i)(j-h)\max_{u\in[i,j]}f''(u)$ by the linear interpolation
remainder. Since $f''(u)=y^2e^{-uy}\le t^{i_0}/(i-i_0)^2$, averaging over
$t=X$ or $t=p_\circ$ gives~\eqref{eq:moment-secant} (the gap is zero at $t=0$).

At $H=256$ the tail satisfies
\begin{equation}\label{eq:power-tail}
 \sum_{h>H}c_h^m\le\tfrac12 c_H^{m-1}
       \left(B_{\rm inv}-1-\sum_{h=1}^H c_h\right).
\end{equation}
Cauchy--Schwarz and monotonicity give
$2c_{H+j}^2\le c_H(c_{H+2j-1}+c_{H+2j})$ for $j\ge1$.
Sum this inequality and use $c_{H+j}^{m-2}\le c_H^{m-2}$.
Equation~\eqref{eq:inverse-identities} then bounds $I(T_r x)$.

In \eqref{eq:fourier-moments} take $n_{\max}=18$, $\Xi=64$ and mesh $1/8$.
Through frequency $8$ use
$\Re\E(1-p_\circ)^{\alpha}-1=\sum_{h\ge1}\Re a_h(\alpha)c_h^m$.
Above order $H$ the coefficients have modulus at most $|a_H(\alpha)|$ since
$a_{h+1}=(h-\alpha)a_h/(h+1)$ and $|h-\alpha|\le h+1$ here for $h\ge H$.
Above frequency $8$ use $0\le\Psi_n(\xi)\le\Psi_n(0)$.
The mode, mesh and frequency errors are bounded as in
Section~\ref{app:insertions}. For the mesh error each step checks
$\E\cosh Z=C_Z=\exp(\lambda\sum_{h\ge2}c_h^m)<5$. The high-frequency error after
$64$ is below $10^{-13}$ for all retained moments.

At a node with $\xi\le8$, let $b$ be the midpoint of the interval for
$kr$ and put $\delta=kr-b$. Let $v$ be the coefficient sum
at the midpoints of the stored $c_h^m$ and $v+d$ the actual coefficient sum.
The exponent of $\Psi_n(\xi)$ is $(b+\delta)(v+d)$. Equation~\eqref{eq:moment-secant} accounts for
interpolation in $d$, and \eqref{eq:power-tail} for its tail.
If $bv\le0$ and $|w|\le d_0<2/5$ for $w=bd+\delta v+\delta d$, enclose
$e^{bv+w}$ by $e^{bv}(1+w)$ with error at most
$\tfrac12 d_0^2\exp(\min(0,bv+d_0))$.
Both $bv$ and the true exponent $bv+w$ are nonpositive, which bounds
the second derivative throughout the Taylor remainder. For each output moment, we sum the coefficients of each shared deviation
($\delta$ or a deviation in $c_1^m,\ldots,c_{24}^m$) over the nodes with
$\xi\le8$. Taking absolute values after summing retains cancellation across nodes. Bound the other linear terms
and the mixed term by absolute values. Use the direct exponential
interval at the remaining nodes.

For invariance at $r\in J\cap H_k$, apply this enclosure to the full
box using the outer panel $J$ for the density bounds. The output lies
inside every upper and positive lower face by
more than $6\cdot10^{-7}$ for $k=3$ and $1.9\cdot10^{-6}$ for $k=4$.
Nonnegativity supplies the zero lower bounds. For entry, start at
$c_h=2^{-h}$ with the same inverse cap. The entry margins exceed
$10^{-8}$ on every panel by the steps stated in the lemma.
The $56$ full-box inclusions and every orbit step include the inverse-cap
check. They use outward rounding on the $2^{-384}$ grid for Fourier
kernels and the $2^{-60}$ grid for orbit products.
\end{proof}

\subsection{Soundness of the finite spectral tests}\label{sec:finite-spectral}

The following $22$ tests bound the full moment operator for every survey
law in its class and every density in the region stated in
Proposition~\ref{prop:finite-contraction}. Only coordinates
with positive diagonal excess contribute to the finite block.

\begin{table}[htp]\centering
\caption{Parameters and strict pivot margins of the finite spectral certificates.}
\label{tab:finite-spectral}
\begin{tabular}{rrrrrrr}\toprule
$k$ & Density tests & $q$ & $\tau$ & $\Xi$ & Block size & Least pivot exceeds\\\midrule
5 & 2  & $0.99$ & $0.02$ & 8  & 14 & $0.023$\\
4 & 16 & $0.99$ & $0.1$  & 8  & 6--17  & $0.035$\\
3 & 4  & $0.98$ & $0.5$  & 12 & 12--15 & $0.044$\\\bottomrule
\end{tabular}
\end{table}

\begin{proof}[Proof of Proposition~\ref{prop:finite-contraction}]
Fix $H=256$. For $k=3$ group the panels in tens, taking the smallest
interval containing each stored moment's bounds on all ten panels.
For $k=4$ use each panel's box. Extend these intervals to all orders
through $H$ by~\eqref{eq:moment-secant}, intersecting with
$0\le c_1\le1/2$ and $c_h\le(B_{\rm inv}-1)/(2h)$ for $h\ge2$.
Let $l_h,u_h$ be the rational endpoints obtained by outward interval
evaluation. For $h>H$ set $u_h=(B_{\rm inv}-1)/(2h)$.
For $k=5$, take outward rational enclosures of the lower bounds
$\mu^h$ from Jensen's inequality and the upper bounds
in~\eqref{eq:invariant-caps} through order $H$, then set
$u_h=7/(8h)$ for $h>H$. In each case use these $u_h$, including the
rounded finite head, in the metric~\eqref{eq:moment-metric}.
They are positive and $O(h^{-1})$, so the metric weights are positive
and summable.

Let $[a,b]$ be the density interval and $q$ the claimed contraction
factor. For $k=3,4$ take $\beta =mkb/2$ and bound $\Psi_n(\xi)$ at $a$,
since it decreases with density for fixed input. Set $\sigma=q$.
For $k=5$ take
$[a,b]=[20.18,22.19]\supset H_5$, put
$\sigma=q/\kappa$, $\kappa=(a+b)^2/(4ab)$, and test the two endpoints
with their own $\beta =mkr/2$ and $\Psi_n(\xi)$. For fixed input, $\Psi_n(\xi)$ is exponential in density, so
$\mathsf K(r)\preceq(1-\theta)\mathsf K(a)+\theta\mathsf K(b)$
at $r=(1-\theta)a+\theta b$ by~\eqref{eq:three-gram}.
Conjugate by $W_0^{1/2}$ and multiply by $mkr/2$. The two endpoint
norms, each with its own density factor, are at most $\sigma$. Thus
the interior norm is at most $\sigma r[(1-\theta)/a+\theta/b]$.
The multiplier equals $1+\theta(1-\theta)(b-a)^2/(ab)\le\kappa$.

For the tabulated $0<\tau<\sigma$, put $d_h=(\beta u_h^{m-1}-\tau)_+$,
$\Lambda=\operatorname{diag}(d_h)$ and $f(y)=(\sigma-\tau)y/(\sigma-\tau y)$.
Define $(\mathcal J v)(n,\xi)=\sum_h f_{n,h}(\xi)v_h$ in measure $g_n(\xi)\,d\xi$ on
each mode, and let $\Omega$ multiply by $\Psi_n(\xi)$. Then
$\mathsf K=\mathcal J^*\Omega \mathcal J$ and the zero-intensity identity gives $\mathcal J^*\mathcal J=I$.
Since $\mathcal J$ is an isometry, $\mathcal J\mathcal J^*$ is an orthogonal projection and
$\mathcal J\mathcal J^*\preceq I$. Together with $\beta W_0\preceq\tau I+\Lambda$, this shows it
suffices that
\begin{equation}\label{eq:spectral-test}
 \Lambda^{1/2}\mathcal J^* f(\Omega)\mathcal J \Lambda^{1/2}\prec(\sigma-\tau) I.
\end{equation}
Indeed, bound $\sqrt\Omega \mathcal J\beta W_0\mathcal J^*\sqrt\Omega$ by
$\tau\Omega+\sqrt\Omega \mathcal J\Lambda \mathcal J^*\sqrt\Omega$. Conjugating its second
term by $(\sigma I-\tau\Omega)^{-1/2}$ and using
$\|B^*B\|=\|BB^*\|$ gives the test. Only $d_h>0$ contribute.

For $n=0,\ldots,8$, let $s_n$ bound $f(\Psi_n(0))$ from above.
Since $f$ is increasing and $\Psi_n(\xi)\le\E {\mathcal V}^n=\Psi_n(0)$ with ${\mathcal V}\le1$,
$s_n$ bounds the whole mode and $s_8$ bounds all modes $n\ge8$.
Normalize each active coordinate $h$ by dividing by $\sqrt{h d_h}$.
Its normalized feature is $\phi_h(\xi)=\Re a_h(n/2+i\xi)$ in measure
$\omega_n(\xi)\,d\xi=g_n(\xi)\,d\xi$. For $n=0$ use instead $\phi_h=\Re a_h(i\xi)/\xi$,
$\omega_0=\xi^2g_0=2g_2$. Put $\mathsf H_n=\int_0^\infty\omega_n\phi\phi^{\mathsf T}\,d\xi$.
Here $\phi=(\phi_h)_h$ in the given mode. The identity $\mathcal J^*\mathcal J=I$ now reads $\sum_{n\ge0}\mathsf H_n=\operatorname{diag}(1/h)$.

Partition the frequency interval $[0,\Xi]$ into cells $I$ of width $1/64$ and set
$\bar\phi_I=|I|^{-1}\int_I\phi(\xi)\,d\xi$. For
$0\le C_I\le\inf_I[(s_n-f(\Psi_n))\omega_n]$, the correction on that
cell is at least $C_I|I|\bar\phi_I\bar\phi_I^{\mathsf T}$ by
Cauchy--Schwarz. Thus the normalized right side minus the left side of
\eqref{eq:spectral-test} is bounded below by
\[
 \operatorname{diag}\!\left(\frac{\sigma-\tau}{h d_h}-\frac{s_8}{h}\right)
 -\sum_{n=0}^7\left[(s_n-s_8)\mathsf H_n-
           \sum_I C_I|I|\bar\phi_I\bar\phi_I^{\mathsf T}\right].
\]
It suffices to show that this lower-bound matrix is positive definite.
We do this by symmetric elimination. A positive lower bound for every pivot
certifies positive definiteness, which explains the last column of Table~\ref{tab:finite-spectral}.

Each $\phi_h$ is a polynomial in $\xi$ with rational coefficients. For $n\ge1$ their full Gram integrals
use
\[
 \int_0^\infty \xi^{2j}g_n(\xi)\,d\xi=\frac{2E_{n,j}}{n\,2^{n+2j}},\qquad
 E_{n,0}=1,\quad E_{n,j+1}=n(n+1)E_{n+2,j}-n^2 E_{n,j}.
\]
This follows by differentiating the cosine transform using
$(\operatorname{sech}^n x)''=n^2\operatorname{sech}^n x-
n(n+1)\operatorname{sech}^{n+2}x$.

To bound $C_I$, enclose $\chi_n(\xi)=\Re\E(1-p_\circ)^{n/2+i\xi}$ at both endpoints
by the binomial expansion through $H$, choosing $l_h^m$ or $u_h^m$
by the coefficient sign. The absolute tail is at most
$|a_H|\sum_{h>H}u_h^m$. For $k=3,4$ sharpen the moment-tail sum
using \eqref{eq:power-tail} with lower bounds through $H$.
For $k=5$ use $(7/8)^4/(3H^3)$. The modulus of the coefficients decreases
after $H$ on these cells. Between endpoints the larger upper bound
increases by at most $|I|^2\E Y^2/8$ since $|\chi_n''|\le\E Y^2$.
Here $\E Y^2=\sum_{h\ge2}2\mathfrak h_{h-1}c_h^m/h$. The same tail
bound applies since $2\mathfrak h_{h-1}\le h$. The gamma recurrence
in Section~\ref{app:insertions} bounds $\omega_n$ from below on each
cell, completing the estimate of $C_I$.

For each of the $22$ tests, outward interval evaluation of these
recurrences and symmetric elimination on the $2^{-384}$ grid gives
the positive pivot bounds in Table~\ref{tab:finite-spectral}.
Thus~\eqref{eq:spectral-test} holds, proving the claimed contraction
factor throughout the stated density region.
\end{proof}

\subsection{Evaluating the endpoint certificates}\label{sec:finite-endpoint-evaluation}

We bound the trial score and squared update residual at $r_-$ and
the trial score at $r_+$, using the same metric as the contraction.

At any evaluation density $t$ we enumerate the possible warning probabilities
$p=\surv_1\cdots\surv_{k-1}$, weighting each tuple of independent surveys
by its probability under the trial survey law. With $Y=-\ln(1-p)$ this gives
$\E Y$, $\E Y^2$ and $\Psi_n(\xi)$ in~\eqref{eq:fourier-site}--\eqref{eq:fourier-moments}.
The clause term uses $k$ draws. The frequency cutoff is $64$ throughout.
At $r_-$ the Fourier mode and moment order cutoffs $(n_{\max},h_{\max})$
are $(32,8192)$ for $k=3$, $(24,256)$ for $k=4$, and $(24,64)$
for $5\le k\le10$. At each $r_+$ take $n_{\max}=32$. Each trial satisfies
\begin{equation}\label{eq:trial-field-cap}
 \lambda\,\E\!\left[\frac{p^2}{1-p}\right]\le\frac85,
 \qquad \lambda=kt/2 .
\end{equation}
Independence of the sign fields gives
$\E\cosh Z=\exp(\lambda\E[p^2/(1-p)])\le e^{8/5}<5$, so $C_Z=5$
suffices in the Fourier error bounds of Section~\ref{app:insertions}.

Write $T=T_t$. For $h\le h_{\max}$, enclose $(Tx)_h$ in
$[\underline m_h,\overline m_h]$ and set $d_h=\max(\overline m_h-x_h,x_h-\underline m_h)$.
For $h>h_{\max}$, monotonicity of both moment sequences gives
$|(Tx)_h-x_h|\le v=\max(\overline m_{h_{\max}},x_{h_{\max}+1})$. For the metric used in the corresponding certificate,
$u_h\le C/h$ in the tail, with $C=11/6$ for $k=3$, $2$ for $k=4$, and
$1$ for $5\le k\le10$. Choose $\bar u_h\ge u_h$ through order $h_{\max}$.
Bounding $\sum_{h>h_{\max}}h^{-(k-1)}$ by its integral gives $\mathcal D_k(Tx,x)^2\le\varepsilon$, where
\[
 \varepsilon=\sum_{h=1}^{h_{\max}} \bar u_h^{\,k-2}d_h^2/h+
       \frac{v^2}{k-2}\left(\frac C{h_{\max}}\right)^{k-2}.
\]
For $k=3,4,5$ and $h\le\min(h_{\max},256)$, take $\bar u_h$ from the checked caps in
Section~\ref{sec:finite-spectral}, using panels 31--40 for $k=3$,
panel 7 for $k=4$, and the common class for $k=5$. Above order 256 at $k=3$
use $11/(6h)$. For $6\le k\le10$ use $\bar u_h=b_h$ through order
$h_{\max}$. These head bounds and the tail constant $C=1$ dominate
$\widehat u_h=\min\{b_h,1/h\}$ and hence bound the residual in $\widehat D_k$.

For $k\ge4$, each operation rounds outwards on the $2^{-384}$ grid.
After range reduction, use the $\operatorname{atanh}$ series for logarithms
and Taylor bounds for exponentials, sines and cosines. Evaluate
$(h+1)a_{h+1}=(h-\alpha)a_h$ at $\alpha=n/2+i\xi$ by interval arithmetic,
bounding the omitted modes separately. For $h_{\max}\le256$ the combined mesh and
high-frequency errors from Section~\ref{app:insertions} are below $10^{-15}$
per moment and $10^{-18}$ for the score.

\begin{proof}[Proof of Lemma~\ref{lem:table-check}]
The class-membership checks were given in
Sections~\ref{sec:finite-scalar-bases} and~\ref{sec:finite-small}.
Substitute the stated atoms and masses in~\eqref{eq:trial-field-cap}
and in the retained Fourier and residual sums above. For the spectral
cases, use the caps defining the corresponding checked metric; for
$k\ge6$, use the head and tail bounds for $\widehat D_k$.

The $k=3$ left endpoint needs the largest moment cutoff.
At $r_-$, enclose $(T_{r_-}x)_h$ for every $h\le8192$ on the frequency
mesh $1/16$. In the binomial recurrence the multiplier at step $h$
has modulus at most $\sqrt{1+(66/h)^2}$. The product over any block of
steps is at most $\exp(33\pi)<2^{160}$, by integrating
$\tfrac12\ln(1+(66/s)^2)$ over $s>0$. Rounding each component to the
$2^{-320}$ grid therefore costs at most $2h/2^{160}$ in each coefficient.
For each moment, the allowance for discarded terms is $10^{-14}$ at each
mode--frequency pair $(n,j/16)$, $0\le n\le32$, $0\le j\le1024$.
Use the Cauchy bound of Section~\ref{app:insertions} for that pair,
or stop the recurrence once the coefficient moduli decrease and the
remaining contributions are below this allowance.

At each $r_-$, add the three Fourier error bounds of Section~\ref{app:insertions}
to the retained score to obtain $u_-$, and insert the moment enclosures
in $\varepsilon$ above. The corrected bound is strictly negative in every row.
At each $r_+$, the omitted modes of the variable cost are nonnegative.
Subtracting $10^{-18}$ for mesh and frequency error leaves $l_+>0$.
At $k=3$ this lower bound exceeds $2.55\cdot10^{-13}$; the corrected
left score is at most $-2.0\cdot10^{-11}$. Rounding down gives
Table~\ref{tab:endpoint-margins}.
\end{proof}

\subsection{Completion of the SP theorem and the enclosures}
\label{sec:five-example}\label{app:proof-endpoint-signs}

\begin{proof}[Proof of Theorem~\ref{thm:uniform-sp} and Proposition~\ref{prop:enclosures}]
For $5\le k\le10$, use the common class on $H_k$; for $k=3,4$, use
the panel classes restricted to $H_k$.
Lemmas~\ref{lem:invariance} and~\ref{lem:entry}, including the finite
scalar checks, give the required compactness, invariance, positive mean
and entry of the initialized orbit.
Lemma~\ref{lem:common-contraction} and
Proposition~\ref{prop:finite-contraction} give contraction in the declared
metrics. Lemma~\ref{lem:table-check} supplies opposite endpoint signs
via the trial-score and residual bounds.
Theorem~\ref{thm:finite-cover-sp} therefore gives the selected nonzero
fixed laws throughout $H_k$ and a unique zero strictly between the
certificate densities in Table~\ref{tab:endpoint-margins}.
The slope bounds and the inward shifts of
Lemma~\ref{lem:enclosure-tightening} then give the sharper enclosures
in Table~\ref{tab:trial-results}.

For $k\ge10$, Proposition~\ref{prop:uniform-obligations} verifies the
same hypotheses, with the uniform bracket of
Lemma~\ref{lem:uniform-endpoints}:
\[
 t_k-\tfrac23 k^22^{-k}<\rstar<t_k-\tfrac1{10}k^22^{-k},
 \qquad t_k=2^k\ln2-(1+\ln2)/2.
\]
Together these cases prove both assertions. The SAT conclusion follows
from the DSS reduction in Section~\ref{sec:update}.
\end{proof}

\section{Evaluation identities and quantitative error bounds}\label{sec:fourier}\label{app:insertions}
\subsection{The Fourier and Gram identities}
\begin{proof}[Proof of Lemma~\ref{lem:fourier-identities}]
Write $C(u,v)=(u+v)/2-\ln(2\cosh((u-v)/2))+\sum_{n\ge1}(e^u+e^v)^{-n}/n$. For $n\ge1$ the beta integral gives
\[
 \int_0^\infty g_n(\xi)\cos(\xi y)\,d\xi=\frac{2}{n[2\cosh(y/2)]^n}.
\]
Thus the $n$th summand is $\tfrac12\int_0^\infty g_n(\xi)\Re e^{-\alpha u-\bar{\alpha}v}\,d\xi$.
Taking expectations introduces $\Psi_n$. For the logarithm,
\[
\frac12\int_0^\infty g_0(\xi)[\cos(\xi(u-v))-1]\,d\xi=-\ln\cosh((u-v)/2),
\]
by differentiating twice and using $\xi^2g_0=2g_2$.

On $e^{-\alpha u-\bar{\alpha}v}$ the multiplier of $\Delta_p$ is $(1-p)^{\alpha}+(1-p)^{\bar{\alpha}}-2=2[\Re(1-p)^{\alpha}-1]$.
Applying $\Delta_p$ to the affine term $(u+v)/2$ gives the constant
$-\ln(1-p)$, which disappears upon a second insertion. In $\E\Delta_p C$ the coefficient of $p^\ell$ is
$2(T_r x)_\ell/\ell$. By the same sign symmetry, in $\E\Delta_p\Delta_y C$
the coefficient of $p^hy^\ell$ is $2B_{\ell h}/\ell$. Equating these coefficients
with those in the Fourier expansion gives
\eqref{eq:fourier-moments} and \eqref{eq:three-gram}.
For $\alpha=n/2+i\xi$ and $0\le p\le1/2$,
\[
 \frac{d^j}{dp^j}(1-p)^{\alpha}
 =(-1)^j \alpha(\alpha-1)\cdots(\alpha-j+1)(1-p)^{\alpha-j}.
\]
Since $\lvert(1-p)^{\alpha-j}\rvert\le 2^j$, the derivative is
$O_j((1+n+\xi)^j)$ uniformly in $p$. Applying the same bound to the factor
involving $y$ gives polynomial growth for each fixed mixed derivative.
For $n\ge1$, the gamma recurrence gives $\xi^2g_n(\xi)\le(n+1)(n+2)g_{n+2}(\xi)$,
while the cosine transform at $y=0$ gives $\int_0^\infty g_n(\xi)\,d\xi=2^{1-n}/n$.
Iterating the inequality $j$ times and then integrating yields
\[
 \int_0^\infty \xi^{2j}g_n(\xi)\,d\xi
 \le \frac{(n+2j)!}{n!}\frac{2^{1-n-2j}}{n+2j}
 =O_j\!\left(\frac{2^{-n}(1+n)^{2j}}{n}\right).
\]
\par\Needspace{9\baselineskip}
\noindent At $n=0$, the insertion factor is $2[\cos(\xi\ln(1-p))-1]$.
For each fixed $j\ge0$, its $j$th derivative in $p$ is even in $\xi$,
vanishes at $\xi=0$, and has bounded second $\xi$-derivative on
$0\le p\le1/2$, $|\xi|\le1$. Taylor's theorem therefore gives
$O_j(\xi^2)$ uniformly in $p$. At zero intensity, inserting a warning of probability $p$ gives root values $p$ and $0$
for the two choices of sign. Thus $B_{\ell h}$ is the coefficient of $p^h$ in $p^\ell$,
so $\mathsf K=I$. For any finitely supported vector $v$, passing from zero
intensity to the given intensity multiplies each nonnegative integrand
$g_n(\xi)(v^{\mathsf T}f_n(\xi))^2$ by $\Psi_n(\xi)\in[0,1]$.
Hence $0\le v^{\mathsf T}\mathsf K v\le\|v\|_2^2$.

\end{proof}

\subsection{Truncation and discretization}
In the variable-cost and moment formulas, retain modes $0\le n\le n_{\max}$.
Replace each integral by a sum on the mesh $\xi=j/16$ up to a frequency cutoff
$\Xi$, with weight $1/16$ at positive nodes and $1/32$ at zero.
For the stated mesh and cutoffs, summing a column of the table bounds
the absolute error in the corresponding quantity. Here $C_Z$ can be any
upper bound on $\E\cosh Z$.
\begin{center}
\begin{tabular}{@{}lcc@{}}
\toprule
Source of error & Updated moment $(T_r x)_\ell$ & Expected variable cost $\E C$\\
\midrule
Modes $n>n_{\max}$ & $4\cdot 2^{-n_{\max}/2}\Psi_{n_{\max}}(0)$ & $2^{-n_{\max}}\Psi_{n_{\max}}(0)/(n_{\max}+1)$\\[2pt]
Mesh $1/16$ & $\dfrac{4(\ell/2+7)C_Z}{e^{16\pi}-1}$ & $\dfrac{10C_Z}{e^{16\pi}-1}$\\[7pt]
Frequencies $\xi>\Xi$ & $\dfrac{9\ell(n_{\max}+1)(3\Xi/2)^{n_{\max}} e^{-\pi \Xi/2}}{\Xi\,n_{\max}!}$ &
$\dfrac{9(n_{\max}+1)(3\Xi/2)^{n_{\max}} e^{-\pi \Xi}}{\Xi\,n_{\max}!}$\\
\bottomrule
\end{tabular}
\end{center}

Admissibility suffices at zero frequency: for independent input surveys,
$Y\le\min(\Lfun(X_1),\Lfun(X_2))$ since $k-1\ge2$.
Writing $q(t)=\Prb(\Lfun(X_1)>t)$ and using $tq(t)\le\E\Lfun(X_1)$ gives
$\E Y^2\le2\int_0^\infty tq(t)^2\,dt\le2(\E\Lfun(X_1))^2<\infty$.

At $\xi=0$, use the limits $g_0(\Psi_0-1)\to-2\lambda\E Y^2/\pi$ and $g_0\Re a_\ell(i\xi)\to-2\mathfrak h_{\ell-1}/(\pi \ell)$, where $\mathfrak h_j=\sum_{i=1}^j1/i$ and $\mathfrak h_0=0$.

At fixed $Z=z$ write $s=(1+e^{-z})^{-1}$, $b=(2\cosh(z/2))^{-1}$ and $f_s(x)=(s-x)/(1-x)$. The sign-averaged moment minus $1/2$ is $[f_s({\mathcal V}b)^\ell+f_{1-s}({\mathcal V}b)^\ell-1]/2$. Write its coefficient of ${\mathcal V}^n$ as $H_{\ell n}(z)$. The first insertion above, before averaging the fields, gives $\widehat H_{\ell n}/\pi=(\ell/2)g_n\Re a_\ell$. On $|x|=(1+s)/2$, $|f_s(x)|\le1$. The radius ratio $2{\mathcal V}b/(1+s)$ is at most ${\mathcal V}/\sqrt2$ for both choices
$s,1-s$. Cauchy's estimate and a geometric sum give the first-row moment bound, using $\E {\mathcal V}^{n_{\max}}=\Psi_{n_{\max}}(0)$. For the variable cost, the omitted terms sum to at most $\sum_{n>n_{\max}}{\mathcal V}^n2^{-n}/n\le {\mathcal V}^{n_{\max}}2^{-n_{\max}}/(n_{\max}+1)$ before averaging. The positive-order terms together have absolute sum at most $2b/(1-1/\sqrt2)$, and $|H_{\ell 0}|\le(\ell/2)e^{-|z|}$. Every truncated sum is thus bounded by $(\ell/2+7)e^{-|z|/2}$.

By Poisson summation the frequency mesh adds copies at $Z+32\pi j$, $j\ne0$. Since $e^{|Z|/2}\le2\cosh Z$, summing the copies gives the second-row moment bound. For the variable cost the positive modes sum to at most $2e^{-|Z|/2}$. For mode zero put $H(z)=-\ln\cosh(z/2)$. By symmetry take $z\ge0$. For every real $s$, write
$\ln(2\cosh(s/2))=s/2+\varepsilon(s)$, where $0\le \varepsilon(s)=\ln(1+e^{-s})\le e^{-s}$.
In the symmetric difference the affine terms cancel. Hence
$\big|[H(z+Z)+H(z-Z)]/2-H(z)\big|$ is at most
$\tfrac12(e^{-(z+Z)}+e^{-(z-Z)})+e^{-z}=e^{-z}(1+\cosh Z)$. Summing over the copies of the zero and positive modes gives the second-row bound for the variable cost.

For the high-frequency tail, Cauchy's estimate on $|p|=1$ gives $|a_\ell|\le2^{n/2}e^{\pi \xi/2}$.

Starting with $g_0=2/[\xi\sinh(\pi \xi)]$, $g_1=2/\cosh(\pi \xi)$, multiply by $(n^2/4+\xi^2)/[(n+1)(n+2)]$ to obtain $g_{n+2}$ from $g_n$. For the cutoffs used here $2(1+n_{\max}^2/(4\Xi^2))\le9/4$, giving $g_n|a_\ell|\le9(3\xi/2)^n e^{-\pi \xi/2}/(\xi\,n!)$ for $\xi\ge \Xi$, $n\le n_{\max}$. The logarithmic derivative of $\xi^{n_{\max}-1}e^{-\pi \xi/2}$ is below $-1/2$. Summing the decreasing tail above $\Xi$ gives the third-row moment bound. For the variable cost the binomial factor is absent, so the exponential decay is $e^{-\pi \xi}$. For a general mesh $1/D$ replace $32\pi$ by $2\pi D$ in the translated copies and $16\pi$ by $\pi D$ in their sum.

\section{Certified refinement}\label{sec:computation}\label{sec:refine}
We count bit operations on a deterministic multitape Turing machine.
We prove Theorem~\ref{thm:precision}.
For $k<16b$ the algorithm repeatedly halves the initial bracket using
a trial score and a slope bound; for $k\ge16b$ the asymptotic enclosure
already gives the required precision. The following approximation
lemma supplies the update used in the first case.
\subsection{Approximating one SP step}\label{app:field-precision}
For $k=3,4$, let $K$ be the class of panel $34$ or $7$, respectively;
for $k\ge5$, use the common class on $H_k$. Use its declared contraction
metric $D$, with the checked common weights for $k\ge6$.
At a rational density $t$ in this panel and precision $p_{\rm dec}\ge2$,
put $P=k+p_{\rm dec}$, $J=\max\{64,P\}$ and $d_{\rm cmp}=4096J+1$.

The numerator and denominator of $t$ have at most $C_0P$ bits for
a fixed encoding constant $C_0$, as in the refinement below.
The input law $x\in K$ has at most $(24P+1)(d_{\rm cmp}+1)+13$ atoms.
Each atom lies on the grid $2^{-(72P+1)}\mathbb Z$ or is a multiple
of $1/1000$, and satisfies $1-v\ge2^{-36P}$, where $v$ is its location. Every mass is a multiple
of $1/\mathfrak d_P$, where $\mathfrak d_P=10^{13}2^{512J+32P}$.

\begin{lemma}[Approximation of one SP step]\label{lem:quantitative-oracle}
Under the preceding assumptions on $K$ and $t$, for integers $k\ge3$
and $p_{\rm dec}\ge2$, from an input $x\in K$ with
the representation just specified one can compute a law
$\widehat x\in K$ with the same representation and
$D(\widehat x,T_tx)<10^{-p_{\rm dec}}$, and an interval of width at most
$10^{-p_{\rm dec}}$ containing $\Gamma_t(x)$, in
$O_{C_0}(P^8\ln^2 P)$ bit operations.
\end{lemma}
\begin{proof}[Construction]

Use the log survey $y=-\ln(1-\surv)$ and the two log fields
$S=-\ln U$, $\bar S=-\ln\overline U$.
For each sign choose the least $s\ge0$ such that $\lambda_0=2^{-s}kt/2\le1$.
Represent a partial warning product $p$ by $z=-\ln(1-p)$.
Multiplication by a survey $c$ sends $z$ to
$-\ln(1-c+ce^{-z})$. Build the warning-product law by successive binary
products, rounding and compressing in this logarithmic coordinate
after each product. A clause with warning probability $p$ contributes $-\ln(1-p)$ to its
log field. Start with a sum of a $\Pois(\lambda_0)$ number of independent
such contributions, transferring the probability of counts above the cutoff
to count zero.
Convolving this law with itself corresponds to adding two independent
such sums. Without approximation, $s$ doublings would recover the
compound Poisson field of mean count $kt/2$. At every addition and doubling, replace each field value $v$ by
$\min(v,1024J)$, retaining its mass.
After each operation, round each support point down, carrying its mass
with it. Within each unit interval keep at most $d_{\rm cmp}+1$ support points
while preserving the moments through degree $d_{\rm cmp}$. To do this, take $d_{\rm cmp}+2$
distinct points $z_i$ with positive weights and subtract from their weights
the largest multiple of $(1/\prod_{j\ne i}(z_i-z_j))_i$ that keeps them
nonnegative. At least one weight becomes zero. Polynomial interpolation
shows that the change integrates every polynomial of degree at most $d_{\rm cmp}$
to zero. Repeat as needed, then round the remaining masses down to the mass
grid and assign the unallocated mass to the point zero. We use the following parameters.
\begin{center}
\begin{tabular}{@{}ll@{}}
\toprule
Poisson count cutoff & $128J$\\
Log-field cutoff & $1024J$\\
Moment-matching degree $d_{\rm cmp}$ & $4096J+1$\\
Log-coordinate and mass grid spacing & $2^{-512J}$\\
Fractional bits for elementary arithmetic & $65536J$\\
Output log-survey cutoff & $24P$\\
Final survey grid spacing & $2^{-(72P+1)}$\\
\bottomrule
\end{tabular}
\end{center}

To bound the accumulated error, change one approximation at a time
and integrate against the function obtained by doing the remaining
operations exactly. For the output moments and score terms, these functions extend
analytically to a Bernstein ellipse of parameter $17/16$ with modulus
at most $2^{40J}$, uniformly in the output moment order. Polynomial
approximation on the ellipse gives error below $2^{-210J}$ per compression,
including rounding. Expanding the doublings into independent copies
gives at most $2^{5J}$ local replacements. The probability of changing
a Poisson count is below $2^{-125J}$ across both fields, and the
positive exponential moment of each field bounds the mass clipped
at any one operation by $2^{-640J}$. Together these estimates give
error below $2^{-120J}$ in every output moment and in each unscaled
score term before final output rounding. Since $(k-1)t\le2^{2J}$,
padding the computed score by $2^{-40P}$ on each side suffices.

After the final clipping, compression and rounding, the error in each
moment is at most $e_0=2^{4-32P}$, and the computed inverse moment
is at most $I(T_t x)+e_0$. Rounding may thus put the law just outside $K$.
We mix it with an explicitly feasible rational law (the anchor) whose
slack in each defining inequality absorbs this error. For $k=3,4$ we use the rational laws supplied for panels 34 and 7,
respectively. For $k\ge5$ we use the point mass at $1/2-2^{-k-2}$. Each satisfies every defining
inequality with slack at least $\sigma=2^{-\max(15,k+2)}$.
The supplement checks the two panel margins with exact rational arithmetic.
All the constraints are affine in the law. After mixing in a fraction
$\theta=e_0/\sigma$ of the anchor, the excess over any upper bound is at
most $(1-\theta)e_0-\theta\sigma\le0$. The same estimate applies to the
shortfall at a lower bound.
The $24P+1$ output cells each contribute at most $d_{\rm cmp}+1$ atoms, with room
for zero and at most 12 anchor atoms. Before mixing, masses are multiples of
$2^{-512J}$. The denominator of $\theta$ divides $2^{32P}$ and the small-panel
anchor masses have denominator $10^{13}$. Hence all output masses lie on the
stated grid. These representation bounds also hold at every larger precision.
The restored law has distance at most
$3(e_0+e_0/\sigma)<10^{-p_{\rm dec}}$ from $T_t x$.

The computation performs $O(P)$ stages. Each produces $O(P^4)$
candidate support points. Per candidate the elimination uses $O(P)$ products,
comparisons and exact divisions on integers of $O(P^2)$ bits. Fast division
as in Brent and Zimmermann~\cite[\S\S1.4.4--1.4.5]{BZ} and fast multiplication
give $O(P^3\ln^2 P)$ bit operations per candidate, hence $O(P^8\ln^2 P)$ bit operations per update and score evaluation.
The detailed error and exact-arithmetic analysis is in the separate
ancillary proof~\cite{SPApproximation}.
\end{proof}

\subsection{Bit complexity of refining the interval}\label{app:proof-refinement}
\begin{proof}[Proof of Theorem~\ref{thm:precision}]
Put $\Delta=10^{-b}$. For $k\ge16b$, Proposition~\ref{prop:enclosures}
gives $|t_k-\rstar|<k^22^{-k}\le2^{-k/2}<\Delta/10$.
Evaluate $t_k=(2^k-1/2)\ln 2-1/2$ at the midpoint of an enclosure of
$\ln 2$ of width $2^{-s}$, $s=k+4b+5$. The resulting $\widehat t$ is
within $\Delta/20$ of $t_k$, so return $(\widehat t-\Delta/2,\widehat t+\Delta/2)$.
The ancillary analysis~\cite{SPApproximation} gives the certified
arithmetic--geometric mean evaluation, using Brent~\cite[\S5]{Brent}
and Brent and Zimmermann~\cite[(4.73)]{BZ}.
Its $O(\ln s)$ steps at $O(s)$ bits cost $O(k\ln^2 k)$ by the multiplication
bound of Harvey and van der Hoeven~\cite{HvH}.

For $k<16b$, start with Table~\ref{tab:trial-results} or, for $k\ge11$,
round~\eqref{eq:asymptotic-window} outwards within $H_k$ to error
$k^22^{-k}/100$ or less, using $\ln 2$ to width $2^{-2k-12}$.
The initial width satisfies $w_0\le1$ and $\ln(1/w_0)=O(k+1)$.
For $k=3,4$ take the invariant class $K$ of the density subinterval that
contains this bracket, namely interval 34 for $k=3$ and interval 7 for $k=4$ (Lemma~\ref{lem:entry}). For $k\ge5$ take the class in Lemma~\ref{lem:invariance}. At initialization and after each halving, write $(l,u)$ for the current
bracket and $w=u-l$ for its width. Before evaluating a midpoint, test whether
$w\le\Delta/10$. If so, return
$(\ell,\ell+\Delta)$ for $\ell=(\Delta/10)\lfloor10l/\Delta\rfloor$,
since $\ell\le l<\rstar<u<\ell+\Delta$.

Let \(D\) be the declared refinement metric on \(K\): \(\widehat D_k\) for \(k\ge6\), and the checked metric of Proposition~\ref{prop:finite-contraction} for \(k=3,4,5\). Write \(a_h\) for its moment caps and \(q\) for its contraction factor.
Couple the Poisson counts at densities $s,t$. For any fixed input, each
output moment changes by at most the probability of an extra clause,
at most $k|s-t|$. The weights in~\eqref{eq:moment-metric} sum to at most
\(1+2\sum_{h\ge2}h^{-2}<3\) since \(a_h\le1\) and \(a_h\le2/h\) for \(h\ge2\)
(Lemma~\ref{lem:contract} and Proposition~\ref{prop:finite-contraction}). Thus $K$ has diameter below $3$. Changing the
density at $\eta_s$ and then applying contraction at $t$ gives
$D(\eta_s,\eta_t)\le3k|s-t|+qD(\eta_s,\eta_t)$, hence
$D(\eta_s,\eta_t)\le3k|s-t|/(1-q)$. We also have $g_k'\ge2^{-k-1}$
by~\eqref{eq:selected-derivative}. For $k=3,4$ use the mean bounds of
Lemma~\ref{lem:entry}. For $k\ge5$ use $(1-2^{-k})^k\ge1-k2^{-k}>1/2$
in~\eqref{eq:invariant-class}.

At the midpoint $t$ maintain a rational $x\in K$ with $D(x,\eta_t)\le aw$,
where $a=2^{-k}/(8k)$. Take the least working precision $p_{\rm dec}\in\mathbb N$ with $p_{\rm dec}\ge2$ and $10^{-p_{\rm dec}}\le aw/1000$.
Section~\ref{app:field-precision} computes an update within $10^{-p_{\rm dec}}$ in $D$
and encloses $\Gamma_t(x)$ in $[\gamma_-,\gamma_+]$ of width at most $10^{-p_{\rm dec}}$.
The maximizing property and~\eqref{eq:interval-comparison} give
\[
 \gamma_-\le g_k(t)\le\gamma_++\tbinom{k}{2}t(aw)^2.
\]
Since $t<2^k$ and $w\le1$, the correction is below $2^{-k}w/128$ and
the score interval has width below $2^{-k-1}w/4$. If its lower bound is
positive, keep the left half of the bracket. If its upper bound is negative,
keep the right half. Otherwise $|g_k(t)|<2^{-k-1}w/4$ and $g_k'\ge2^{-k-1}$
give $|t-\rstar|<w/4$, so keep the central half.

To initialize $x$, start from the feasible rational law (the anchor) in
Section~\ref{app:field-precision}. Perform $j$ approximate updates, taking the least $j$ with
$3(0.99)^j\le aw_0/2$, so $j=O(k+1)$. Errors from approximate updates
sum to at most $10^{-p_{\rm dec}}/(1-q)\le aw/10$, giving initial distance at most
$3aw_0/5$. After a halving to $w'=w/2$ the midpoint moves by at most $w/4$.
By the density estimate and $1-q\ge0.01$, the old trial is within
$(2+1200k^22^k)aw'$ of the new fixed law. At the new precision, 2048
updates suffice because $q^{2048}(2+1200k^22^k)\le1/2$. For $k\le5$ use
$q\le0.99$, and for $k=6$ use $q<0.8$. Thereafter $q_{k+1}/q_k\le0.784$
by~\eqref{eq:common-factor}, and $0.784^{2048}\,2(7/6)^2<1$.
Including rounding gives distance at most $3aw'/5$.

There are $O(b)$ halvings and $O(k+b)=O(b)$ calls. Each has $w>\Delta/10$,
so $p_{\rm dec}=O(k+b)=O(b)$. After $i$ halvings the endpoints are dyadic subdivisions
of the initial rational bracket and use $O(k+i)$ bits, hence $O(k+p_{\rm dec})$
since $w=w_0/2^i\le2^{-i}$. Section~\ref{app:field-precision} gives
$O((k+p_{\rm dec})^8\ln^2(k+p_{\rm dec}))$ bit operations per call, or $O(b^9\ln^2 b)$ in all.
\end{proof}

\end{document}